\documentclass[10pt,reqno]{amsart}

\ifdefined\pdfsuppresswarningpagegroup
\fi
\ifdefined\pdfminorversion
\fi

\usepackage{amsmath,amssymb,amsthm,mathtools}
\usepackage{enumitem}
\usepackage{xcolor}
\usepackage{hyperref}
\usepackage{bookmark}

\theoremstyle{definition}
\newtheorem{theorem}{Theorem}[section]
\newtheorem{lemma}[theorem]{Lemma}
\newtheorem{proposition}[theorem]{Proposition}
\newtheorem{corollary}[theorem]{Corollary}
\newtheorem*{maintheorema}{Theorem A}
\usepackage{longtable,pdflscape,array}
\usepackage[pass]{geometry}
\newtheorem*{maintheoremb}{Theorem B}

\theoremstyle{definition}
\newtheorem{definition}[theorem]{Definition}
\newtheorem{remark}[theorem]{Remark}

\providecommand{\sslash}{\mathbin{/\mkern-6mu/}}

\newcommand{\SO}{\operatorname{SO}}
\newcommand{\Orth}{\operatorname{O}}
\newcommand{\Id}{\mathrm{id}}
\newcommand{\Cyc}[1]{\mathrm C_{#1}}
\newcommand{\Dih}[1]{\mathrm D_{#1}}
\newcommand{\BDih}[1]{\mathrm D_{#1}^{*}}
\newcommand{\familyref}[1]{\textup{Family~\ref{#1}}}
\newcommand{\familiesref}[2]{\textup{Families~\ref{#1} and~\ref{#2}}}

\hypersetup{
colorlinks=true,
linkcolor=blue,
citecolor=blue,
urlcolor=blue,
pdftitle={Orbifold Morse Handles I. Morse Indices and Isomorphism Classes of Ordered Handles},
pdfauthor={Suhyoung Choi and Juseop Lee},
pdfsubject={Invariant splittings determined by Hessians, Morse indices, and isomorphism classes of ordered handles},
pdfkeywords={orbifold Morse theory, orbifold handles, finite subgroups of SO(4), invariant splittings determined by Hessians}
}

\title[Orbifold Morse Handles I]{Orbifold Morse Handles I.\\Morse Indices and Isomorphism Classes of Ordered Handles}
\author{Suhyoung Choi}
\address{Department of Mathematical Sciences, KAIST, 291 Daehak-ro Yuseong-gu Daejeon, South Korea}
\email{schoi@math.kaist.ac.kr}
\author{Juseop Lee}
\address{Department of Mathematical Sciences, KAIST, 291 Daehak-ro Yuseong-gu Daejeon, South Korea}
\email{dlwntjq77@kaist.ac.kr}
\date{}

\begin{document}

\begin{abstract}
We classify the Morse indices and ordered handles arising from nondegenerate critical points of smooth functions on effective orientable four-dimensional orbifolds.
Let \(\rho\colon G\to\SO(4)\) be a faithful representation of a finite group.
Let \(V\subset\mathbb R^4\) be the negative eigenspace of a nondegenerate Hessian that is invariant under \(\rho(G)\).
The splitting \(\mathbb R^4=V\oplus V^\perp\) determines an ordered quotient handle.
Two such ordered handles admit an orientation-preserving isomorphism if and only if their negative eigenspaces lie in the same orbit under \(N_{\SO(4)}(\rho(G))\).
We determine the possible Morse indices and the normalizer orbits for every family in the completed Du Val list.
\end{abstract}

\maketitle
\enlargethispage{3pt}
\tableofcontents

\section{Introduction}\label{sec:introduction}

Classical Morse theory describes the change in sublevel sets across a critical value by a handle attachment~\cite{milnor1963morse}.
Let \(x\) be a critical point of an orbifold Morse function in the sense of Definition~\ref{def:orbifold-morse-critical}.
Let \(G_x\) be its local group, as defined in Section~\ref{subsec:orbifolds}.
The group \(G_x\) acts on the tangent space at a lift of \(x\).
The Hessian is \(G_x\)-invariant.
With respect to a \(G_x\)-invariant inner product, the negative and positive eigenspaces of the Hessian are \(G_x\)-invariant.
By the orbifold Morse lemma, the function has a quadratic normal form with these eigenspaces as its negative and positive subspaces~\cite[Theorem~4.7 and Lemma~4.9]{hepworth2009morse}.
In dimension four, the local classification reduces to invariant orthogonal splittings of \(\mathbb R^4\) under finite subgroups of \(\SO(4)\).

We classify the invariant negative eigenspaces and the ordered handles defined in Definition~\ref{def:ordered-handle}.
The Morse index is the dimension of the negative eigenspace.
For a finite group \(G\) and a faithful representation \(\rho\colon G\to\SO(4)\), the isomorphism classes of ordered handles are in bijection with the orbits of \(N_{\SO(4)}(\rho(G))\) on the invariant negative eigenspaces.
Invariant negative eigenspaces of the same dimension need not lie in the same orbit.
We compute these orbits for every family in the completed Du Val list in Table~\ref{tab:mecchia-seppi-families}.

\subsection{Main results}\label{subsec:main-results}

Throughout this paper, we assume that orbifolds are smooth, effective, and orientable.
If a finite group \(G\) acts smoothly and effectively on a manifold \(X\), we write \(X\sslash G\) for the quotient orbifold.
If \(G\) acts smoothly on manifolds \(X\) and \(Y\) and the diagonal action on \(X\times Y\) is effective, we write
\[
X\times_GY=(X\times Y)\sslash G.
\]
If \(X\times Y=\varnothing\), we set \(X\times_GY=\varnothing\).
For a Euclidean space \(V\), let \(D(V)\) and \(S(V)\) denote its unit disk and sphere.
For an orbifold \(\mathcal X\), write \(|\mathcal X|\) for its underlying space.

\begin{definition}\label{def:ordered-handle}
Let \(G\) be a finite group, and let \(\rho\colon G\to\SO(4)\) be faithful.
For \(0\le k\le4\), let \(\operatorname{Gr}_k(\mathbb R^4)\) denote the Grassmannian of \(k\)-dimensional linear subspaces of \(\mathbb R^4\).
Set
\[
\operatorname{Gr}_{k}(\rho)
=
\{V\in\operatorname{Gr}_{k}(\mathbb R^{4})\mid \rho(G)V=V\}.
\]
Let \(V\in\operatorname{Gr}_{k}(\rho)\) and set \(V^+=V^\perp\).
We use the convention that \(D(0)\) is a point and \(S(0)=\varnothing\).
The quotient handle determined by \(V\) is
\[
H_\rho(V)=D(V)\times_GD(V^+).
\]
Let
\[
q_{\rho,V}\colon D(V)\times D(V^+)\longrightarrow |H_\rho(V)|
\]
denote the quotient map.
The attaching and belt regions of the quotient handle are
\[
A_\rho(V)=S(V)\times_GD(V^+),
\qquad
B_\rho(V)=D(V)\times_GS(V^+).
\]
The core and cocore of the quotient handle are
\[
K_\rho(V)=q_{\rho,V}\bigl(D(V)\times\{0\}\bigr),
\qquad
K^*_\rho(V)=q_{\rho,V}\bigl(\{0\}\times D(V^+)\bigr).
\]
Let \(o_{\rho,V}\) denote the image of the origin.
The \emph{ordered handle} \(\mathbf H_\rho(V)\) consists of \(H_\rho(V)\) together with \(A_\rho(V)\), \(B_\rho(V)\), \(K_\rho(V)\), \(K^*_\rho(V)\), and \(o_{\rho,V}\).
\end{definition}

\begin{definition}\label{def:ordered-handle-isomorphism}
Let \(G\) be a finite group, and let \(\rho\colon G\to\SO(4)\) be faithful.
Let \(V\in\operatorname{Gr}_{k}(\rho)\) and \(V'\in\operatorname{Gr}_{k'}(\rho)\), where \(0\le k,k'\le4\).
An \emph{isomorphism}
\[
\mathbf H_\rho(V)\longrightarrow\mathbf H_\rho(V')
\]
is an orientation-preserving orbifold diffeomorphism \(F\colon H_\rho(V)\to H_\rho(V')\), as in Definition~\ref{def:orbifold-diffeomorphism-convention}, satisfying
\[
\begin{aligned}
F(o_{\rho,V})&=o_{\rho,V'},&
F(A_\rho(V))&=A_\rho(V'),&
F(B_\rho(V))&=B_\rho(V'),\\
F(K_\rho(V))&=K_\rho(V'),&
F(K^*_\rho(V))&=K^*_\rho(V').
\end{aligned}
\]
The \emph{ordered dual} of \(\mathbf H_\rho(V)\) is \(\mathbf H_\rho(V^\perp)\).
\end{definition}

For a faithful representation \(\rho\colon G\to\SO(4)\) of a finite group, set \(N^+_\rho=N_{\SO(4)}(\rho(G))\).
For \(0\le k\le4\), the group \(N^+_\rho\) acts on \(\operatorname{Gr}_k(\rho)\) by \(T\cdot V=T(V)\).
The following theorem identifies these orbits with the classes of ordered handles under orientation-preserving isomorphism.

\begin{maintheorema}
Let \(\rho\colon G\to\SO(4)\) be a faithful representation of a finite group.
For \(0\le k\le4\), the assignment
\[
V\longmapsto\mathbf H_\rho(V)
\]
induces a bijection between the \(N^+_\rho\)-orbits in \(\operatorname{Gr}_k(\rho)\) and the isomorphism classes of ordered handles determined by \(\rho(G)\)-invariant \(k\)-dimensional subspaces.
Each isomorphism class is realized by a \(\rho(G)\)-invariant quadratic function on \(\mathbb R^4\) with a nondegenerate critical point at the origin.
For \(V\in\operatorname{Gr}_k(\rho)\), the ordered dual of the class represented by \(V\) is represented by \(V^\perp\).
\end{maintheorema}

Proposition~\ref{prop:ordered-handle-orbits} identifies the normalizer orbits with the isomorphism classes of ordered handles.
Lemma~\ref{lem:orbifold-morse-normal-form} identifies the invariant splitting at a nondegenerate critical point.
Lemma~\ref{lem:main-hessian-splitting-principle} proves that every invariant orthogonal splitting is realized by an invariant quadratic function.
Theorem A is proved at the end of Section~\ref{sec:local-conventions}.

Theorem~\ref{thm:4d-morse-index-classification} assigns every family in the completed Du Val list, for each admissible choice of parameters, to one of the following four cases.
Every summand in the first three rows is irreducible.
In the last row, the two-dimensional summand may split into two invariant lines.
\[
\begin{array}{c|c}
\text{dimensions of the summands}
&
\text{possible Morse indices}
\\ \hline
4 & \{0,4\}\\
2\oplus2 & \{0,2,4\}\\
1\oplus3 & \{0,1,3,4\}\\
1\oplus1\oplus2 & \{0,1,2,3,4\}
\end{array}.
\]

The representatives of the families and the admissible parameters are fixed in Section~\ref{subsec:so4-finite-subgroups} and Definition~\ref{def:representative-family-parameters}.
For \textup{Families} \(1'\) and \(11'\), write \(r=2r'\).
All spans are over \(\mathbb R\) in the quaternionic basis \((1,\mathbf{i},\mathbf{j},\mathbf{k})\).

For every finite \(G<\SO(4)\), the unique representatives for indices \(0\) and \(4\) are \(0\) and \(\mathbb H\), respectively.

Each subspace \(V\) in the last column of Table~\ref{tab:normalizer-orbit-representatives} represents the splitting \(V\oplus V^+\), where \(V^+=V^\perp\).
For each group, distinct entries represent distinct orbits.
The representatives for index \(3\) are the orthogonal complements of the representatives for index \(1\).
If a group has no entry for index \(1\), it has no invariant line.

The rows for the trivial group are separate from those for the parameterized families.
Under factor interchange, quaternionic conjugation sends the subspaces in the table to representatives of the normalizer orbits for the resulting group, as in Definition~\ref{def:factor-interchange}.

\clearpage
\newgeometry{margin=18mm}
\begin{landscape}
\begingroup
\normalsize
\setlength{\tabcolsep}{6pt}
\renewcommand{\arraystretch}{1.25}
\begin{longtable}{@{}>{\raggedright\arraybackslash}p{0.19\linewidth}>{\centering\arraybackslash}p{0.08\linewidth}>{\raggedright\arraybackslash}p{0.34\linewidth}>{\raggedright\arraybackslash}p{0.32\linewidth}@{}}
\caption{Representatives of the normalizer orbits of invariant negative eigenspaces.}\label{tab:normalizer-orbit-representatives}\\

\hline
Group or family & Dimension \(k\) & Condition & Representatives \(V\)\\
\hline
\endfirsthead
\multicolumn{4}{c}{\tablename~\thetable\ continued}\\
\hline
Group or family & Dimension \(k\) & Condition & Representatives \(V\)\\
\hline
\endhead
\hline
\endfoot
\(G=\{1\}\) & \(1\) & --- & \(\operatorname{span}\{1\}\)\\
\(G=\{1\}\) & \(2\) & --- & \(\operatorname{span}\{1,\mathbf{i}\}\)\\
\hline
\(1,11\) & \(1\) & \(m=n=1,\ r\mid(1-s)\) & \(\operatorname{span}\{1\}\)\\
\(1,11\) & \(1\) & \(m=n=1,\ r\nmid(1-s),\ r\mid(1+s)\) & \(\operatorname{span}\{\mathbf{j}\}\)\\
\hline
\(1'\) & \(1\) & \(m=n=1,\ r'=2\) & \(\operatorname{span}\{1\},\quad\operatorname{span}\{\mathbf{j}\}\)\\
\(1'\) & \(1\) & \(m=n=1,\ r'\ne2,\ r'\mid(1-s)\) & \(\operatorname{span}\{1\}\)\\
\(1'\) & \(1\) & \(m=n=1,\ r'\ne2,\ r'\nmid(1-s),\ r'\mid(1+s)\) & \(\operatorname{span}\{\mathbf{j}\}\)\\
\hline
\(11'\) & \(1\) & \(m=n=1,\ r'=2\) & \(\operatorname{span}\{1\},\quad\operatorname{span}\{\mathbf{j}\}\)\\
\(11'\) & \(1\) & \(m=n=1,\ r'>2,\ (1-s)/r'\in2\mathbb Z\) & \(\operatorname{span}\{1\},\quad\operatorname{span}\{\mathbf{i}\}\)\\
\(11'\) & \(1\) & \(m=n=1,\ r'>2,\ (1-s)/r'\in2\mathbb Z+1\) & \(\operatorname{span}\{1\}\)\\
\(11'\) & \(1\) & \(m=n=1,\ r'>2,\ (1+s)/r'\in2\mathbb Z\) & \(\operatorname{span}\{\mathbf{j}\},\quad\operatorname{span}\{\mathbf{k}\}\)\\
\(11'\) & \(1\) & \(m=n=1,\ r'>2,\ (1+s)/r'\in2\mathbb Z+1\) & \(\operatorname{span}\{\mathbf{j}\}\)\\
\hline
\pagebreak
\(1\) & \(2\) & \(r=1\) & \(\operatorname{span}\{1,\mathbf{i}\}\)\\
\(1\) & \(2\) & \(r=2,\ (m,n)\ne(1,1)\) & \(\operatorname{span}\{1,\mathbf{i}\}\)\\
\(1\) & \(2\) & \(m=n=1,\ r=2\) & \(\operatorname{span}\{1,\mathbf{i}\},\quad\operatorname{span}\{1,\mathbf{j}\}\)\\
\(1\) & \(2\) & \(r>2\) & \(\operatorname{span}\{1,\mathbf{i}\},\quad\operatorname{span}\{\mathbf{j},\mathbf{k}\}\)\\
\hline
\(1'\) & \(2\) & \(r'=1\) & \(\operatorname{span}\{1,\mathbf{i}\}\)\\
\(1'\) & \(2\) & \(m=n=1,\ r'=2\) & \(\operatorname{span}\{1,\mathbf{i}\},\quad\operatorname{span}\{\mathbf{j},\mathbf{k}\},\quad\operatorname{span}\{1,\mathbf{j}\}\)\\
\(1'\) & \(2\) & \(r'>1,\ (m,n,r')\ne(1,1,2)\) & \(\operatorname{span}\{1,\mathbf{i}\},\quad\operatorname{span}\{\mathbf{j},\mathbf{k}\}\)\\
\hline
\(11\) & \(2\) & \(r\in\{1,2\}\) & \(\operatorname{span}\{1,\mathbf{i}\}\)\\
\(11\) & \(2\) & \(r>2\) & \(\operatorname{span}\{1,\mathbf{i}\},\quad\operatorname{span}\{\mathbf{j},\mathbf{k}\}\)\\
\hline
\(11'\) & \(2\) & \(r'=1\) & \(\operatorname{span}\{1,\mathbf{i}\}\)\\
\(11'\) & \(2\) & \(r'>1\) & \(\operatorname{span}\{1,\mathbf{i}\},\quad\operatorname{span}\{\mathbf{j},\mathbf{k}\}\)\\
\hline
\(3\) & \(2\) & \(m=1,\ n\ge2\) & \(\operatorname{span}\{1+\mathbf{j},\mathbf{i}-\mathbf{k}\}\)\\
\(34\) & \(2\) & \(m=1,\ n\ge3,\ n\text{ odd}\) & \(\operatorname{span}\{1+\mathbf{j},\mathbf{i}-\mathbf{k}\}\)\\
\(3\) & \(1,2,3\) & \(m>1,\ n\ge2\) & \(\varnothing\)\\
\(34\) & \(1,2,3\) & \(m>1,\ n\ge3,\ m,n\text{ odd}\) & \(\varnothing\)\\
\hline
\(21,21',26,26',\) \(31,31'\) & \(1\) & --- & \(\operatorname{span}\{1\}\)\\
\(21,21',26,26',\) \(31,31'\) & \(2\) & --- & \(\varnothing\)\\
\(26''\) & \(1\) & --- & \(\operatorname{span}\{\mathbf{i}\}\)\\
\(26''\) & \(2\) & --- & \(\varnothing\)\\
\hline
\(2,4,5\text{--}10,12\text{--}20,\) \(22\text{--}25,27\text{--}30,\) \(32,32',33,33'\) & \(1,2,3\) & All parameters & \(\varnothing\)\\
\hline
\end{longtable}
\endgroup
\end{landscape}
\restoregeometry
\clearpage

\begin{maintheoremb}
Let \(G<\SO(4)\) be one of the representatives fixed in Section~\ref{subsec:so4-finite-subgroups}, and let \(\rho\colon G\hookrightarrow\SO(4)\) be the inclusion.
For \(k=1,2\), Table~\ref{tab:normalizer-orbit-representatives} lists exactly one representative of each \(N^+_\rho\)-orbit in \(\operatorname{Gr}_k(\rho)\).
The representatives in \(\operatorname{Gr}_3(\rho)\) are the orthogonal complements of those in \(\operatorname{Gr}_1(\rho)\).
The unique representatives in dimensions \(0\) and \(4\) are \(0\) and \(\mathbb H\), respectively.
The associated ordered handles \(\mathbf H_\rho(V)\) form a complete set of representatives for the isomorphism classes of ordered handles determined by \(\rho\).
Each class is realized by an invariant quadratic function with a nondegenerate critical point at the origin.
\end{maintheoremb}

Theorem B is proved at the end of Section~\ref{sec:four-dimensional-handle-classification}.

\subsection{Previous work and organization}\label{subsec:relation-previous-work}

Wasserman developed Morse theory for invariant functions and described changes across critical values by handle-bundles over critical orbits~\cite{wasserman1969equivariant}.
Lerman and Tolman proved a Morse lemma for nondegenerate critical suborbifolds.
They also described, up to homotopy, the change across a Bott-Morse critical value by attaching the disk bundle of the negative normal orbi-bundle~\cite[Lemmas~4.3 and~4.4]{lermantolman1997hamiltonian}.
Hepworth defined the index and coindex as isomorphism classes of representations of the automorphism group.
He proved a Morse lemma and an attachment theorem up to homotopy for the underlying space of a differentiable Deligne-Mumford stack~\cite[Definition~4.5 and Theorems~4.7 and~7.6]{hepworth2009morse}.
Cho and Hong constructed Morse-Smale-Witten complexes from Morse-Smale functions on effective orientable orbifolds~\cite{chohong2014orbifold}.
Corrigan established Morse inequalities for integral orbifold Borel homology~\cite{corrigan2020morse}.
Ortiz and Valencia proved that the Morse condition for Lie groupoid morphisms is invariant under Morita equivalence~\cite[Proposition~3.8]{ortizvalencia2024morse}.
They also proved a Morse lemma for proper Lie groupoids~\cite[Theorem~5.2]{ortizvalencia2024morse}.
Their attachment theorem describes the change in sublevel groupoids across a nondegenerate critical orbit~\cite[Theorem~6.7]{ortizvalencia2024morse}.
Under the hypotheses of that theorem, the groupoid obtained by attaching the negative disk groupoid to the lower sublevel groupoid is a deformation retract of the upper sublevel groupoid.

These works do not classify the invariant negative eigenspaces for a fixed oriented isotropy representation.
We classify the invariant negative eigenspaces and the ordered handles determined by these eigenspaces for every finite subgroup of \(\SO(4)\).

Threlfall and Seifert classified the finite subgroups of \(\SO(4)\) up to conjugacy~\cite{seifertthrelfall1931,seifertthrelfall1933}.
Du Val organized these groups into numbered quaternionic families~\cite{val1964homographies}.
Mecchia and Seppi added \textup{Families} \(33\), \(33'\), and \(34\) to Du Val's list.
Table~\ref{tab:mecchia-seppi-families} lists the completed Du Val families.
For \textup{Family} \(33'\), we use the correction to the quotient isomorphism in Remark~\ref{rem:family-33p-angle-correction}.
Mecchia and Seppi determined \(N_{\SO(4)}(G)/G\) by computing normalizers in \(S^3\times S^3\)~\cite[Section~3.1]{mecchiaseppi2019isometry}.
The action of \(N^+_\rho\) on \(\operatorname{Gr}_k(\rho)\) factors through \(N^+_\rho/\rho(G)\).
The abstract isomorphism type of this quotient does not specify the induced action.
We determine the \(N^+_\rho\)-orbits in \(\operatorname{Gr}_k(\rho)\) directly.
Proposition~\ref{prop:ordered-handle-orbits} identifies these orbits with the isomorphism classes of ordered handles.
Rastanawi and Rote characterized the elements of \(\SO(4)\) that preserve a given great circle in \(S^3\)~\cite[Proposition~4.5]{rastanawirote2022pointgroups}.
For representations of a finite group in \(\Orth(4)\), Garc\'ia-Hern\'andez and Williams proved that conjugacy by an orientation-preserving diffeomorphism of \(S^3\) implies conjugacy by an element of \(\SO(4)\)~\cite{garciahernandezwilliams2025}.

Section~\ref{sec:local-conventions} recalls the orbifold conventions, defines quotient handles, and reduces the classification to invariant splittings.
Section~\ref{sec:morse-index-classification} introduces the quaternionic notation for the families and proves Theorem~\ref{thm:4d-morse-index-classification}.
Section~\ref{sec:four-dimensional-handle-classification} determines the normalizer orbits and proves Theorem B.
Appendix~\ref{app:morse-index-case-analysis} contains the calculations with generators for the four-dimensional families.

\section{Preliminaries}\label{sec:local-conventions}

In this section, we relate nondegenerate critical points to invariant splittings and ordered quotient handles.
Section~\ref{subsec:orbifolds} fixes the conventions for orbifold charts, smooth functions, and orbifold diffeomorphisms.
Section~\ref{subsec:hessian-splittings} describes the splittings determined by Hessians and their realization by invariant quadratic functions.
Section~\ref{subsec:main-quotient-handle} identifies the isomorphism classes of ordered handles with the normalizer orbits and proves Theorem A.

\subsection{Orbifold charts and smooth functions}\label{subsec:orbifolds}

We use effective smooth orbifolds in the classical sense~\cite[Chapter~2]{cooper2000three}.
Every point of an \(n\)-orbifold has a neighborhood \(U\) represented by a chart \((\widetilde U,\Gamma,\pi_U)\).
The domain \(\widetilde U\subset\mathbb R^n\) is connected and open.
The finite group \(\Gamma\) acts smoothly and effectively on \(\widetilde U\).
The map \(\pi_U\colon \widetilde U\longrightarrow U\) is a \(\Gamma\)-invariant continuous surjection inducing a homeomorphism \(\widetilde U/\Gamma\cong U\).
After passing to common refinements, changes of charts are represented by equivariant smooth embeddings.

\begin{definition}\label{def:local-group-orientability}
Let \(x\in|\mathcal O|\).
Choose a chart \((\widetilde U,\Gamma,\pi_U)\) around \(x\) and a lift \(\widetilde x\in\widetilde U\) with \(\pi_U(\widetilde x)=x\).
The \emph{local group} \(G_x\) is the stabilizer of \(\widetilde x\) in \(\Gamma\).
Its isomorphism type is independent of the chart and the lift.
An orbifold is \emph{orientable} if its charts admit orientations preserved by every local group and every chart embedding.
\emph{Orbifolds with boundary} and \emph{orbifolds with corners} are defined using half-space and orthant charts, respectively.
The quotient handles in this paper are orbifolds with corners.
Their orbifold diffeomorphisms preserve the corner structure.
\end{definition}

\begin{definition}\label{def:orbifold-diffeomorphism-convention}
Let \(\mathcal O\) and \(\mathcal P\) be smooth orbifolds.
A \emph{smooth orbifold map} \(F\colon\mathcal O\longrightarrow\mathcal P\) is a continuous map of the underlying spaces admitting smooth local lifts in compatible orbifold charts.
The local lifts are equivariant for homomorphisms between the corresponding local groups.
They are compatible after passing to common refinements.

An \emph{orbifold diffeomorphism} is a smooth orbifold map with a smooth orbifold inverse.
If \(\mathcal O\) and \(\mathcal P\) are oriented, an orbifold diffeomorphism is \emph{orientation-preserving} if its local lifts preserve orientation.
\end{definition}

\begin{remark}\label{rem:orthogonal-local-groups}
The tangent representation of \(G_x\) is faithful.
With respect to a \(G_x\)-invariant inner product and an orthonormal basis, it identifies \(G_x\) with a finite subgroup of \(\Orth(n)\)~\cite[Theorem~2.3]{cooper2000three}.
For an orientable orbifold, this subgroup lies in \(\SO(n)\).
In dimension four, we identify \(G_x\) with the corresponding finite subgroup of \(\SO(4)\).
\end{remark}

\begin{definition}\label{def:real-valued-orbifold-function}
Let \(\mathcal O\) be a smooth orbifold.
A \emph{real-valued orbifold function} is a continuous function \(f\colon|\mathcal O|\to\mathbb R\) such that \(f\circ\pi_U\colon\widetilde U\to\mathbb R\) is smooth for every orbifold chart \((\widetilde U,\Gamma,\pi_U)\).
\end{definition}

The local functions \(f\circ\pi_U\) are \(\Gamma\)-invariant.
They agree under changes of charts.
Since \(\mathbb R\) has trivial isotropy, this definition agrees with that of a smooth orbifold map to \(\mathbb R\) in Definition~\ref{def:orbifold-diffeomorphism-convention}.

\subsection{Orbifold Morse functions and invariant splittings}\label{subsec:hessian-splittings}

\begin{definition}\label{def:orbifold-morse-critical}
Let \(\mathcal O\) be a smooth orbifold, and let \(f\colon |\mathcal O|\to\mathbb R\) be a real-valued orbifold function in the sense of Definition~\ref{def:real-valued-orbifold-function}.
Let \(x\in|\mathcal O|\).
Choose an orbifold chart \((\widetilde U,G_x,\pi_U)\) around \(x\) and a lift \(\widetilde x\in\widetilde U\) fixed by \(G_x\).
Set \(\widetilde f=f\circ\pi_U\colon\widetilde U\to\mathbb R\).

The point \(x\) is a \emph{critical point} of \(f\) if \(d\widetilde f_{\widetilde x}=0\).
The function \(f\) is \emph{Morse at \(x\)} if \(x\) is a critical point and the Hessian of \(\widetilde f\) at \(\widetilde x\) is nondegenerate.
If \(f\) is Morse at \(x\), the \emph{Morse index of \(f\) at \(x\)} is the number of negative eigenvalues of this Hessian, counted with multiplicity.
The function \(f\) is an \emph{orbifold Morse function} if it is Morse at every critical point.
\end{definition}

By the chain rule, the condition \(d\widetilde f_{\widetilde x}=0\) is independent of the chart and lift.
At a critical point, a change of chart or lift transforms the Hessian by congruence.
By Sylvester's law of inertia, the Morse condition and the Morse index are independent of the chart and lift.

\begin{proposition}\label{prop:local-morse-general}
Let \(G<\Orth(n)\) be finite, and let \(F\colon\mathbb R^n\to\mathbb R\) be a \(G\)-invariant smooth function with a nondegenerate critical point at the origin.
Then there is a \(G\)-invariant orthogonal decomposition
\[
\mathbb R^n=V_1\oplus\cdots\oplus V_s
\]
into irreducible real \(G\)-submodules such that the Hessian of \(F\) is definite on each \(V_i\).
The Morse index is the sum of the dimensions of the summands on which the Hessian is negative definite.
If the real \(G\)-module \(\mathbb R^n\) is irreducible, then the Morse index is \(0\) or \(n\).
\end{proposition}

\begin{proof}
Let \(\mathsf B\) be the Hessian of \(F\) at the origin.
With respect to a \(G\)-invariant inner product, write \(\mathsf B(v,w)=\langle \mathcal A v,w\rangle\).
The operator \(\mathcal A\) is self-adjoint and invertible.
It commutes with the action of \(G\).
Its eigenspaces are \(G\)-invariant and pairwise orthogonal.
Each eigenspace admits an orthogonal decomposition into irreducible real \(G\)-submodules.
The sign of the Hessian on each summand is the sign of the corresponding eigenvalue.
The Morse index is the sum of the dimensions of the summands on which the Hessian is negative definite.
If \(\mathbb R^n\) is irreducible, then the Hessian is definite and the Morse index is \(0\) or \(n\).
\end{proof}

Fix \(f\), \(x\), \((\widetilde U,G_x,\pi_U)\), \(\widetilde x\), and \(\widetilde f\) as in Definition~\ref{def:orbifold-morse-critical}, and assume that \(f\) is Morse at \(x\).
Choose a \(G_x\)-invariant inner product on \(T_{\widetilde x}\widetilde U\) and identify this Euclidean space with \(\mathbb R^n\).

\begin{lemma}\label{lem:orbifold-morse-normal-form}
The Hessian of \(\widetilde f\) at \(\widetilde x\) determines a \(G_x\)-invariant orthogonal splitting \(\mathbb R^n=V^-\oplus V^+\).
There are \(G_x\)-equivariant local coordinates \((v^-,v^+)\in V^-\oplus V^+\) centered at \(\widetilde x\) in which
\[
\widetilde f(v^-,v^+)
=
f(x)-\lVert v^-\rVert^2+\lVert v^+\rVert^2.
\]
The Morse index of \(f\) at \(x\) is \(\dim V^-\).
\end{lemma}

\begin{proof}
Let \(\mathsf B\) be the Hessian of \(\widetilde f\) at \(\widetilde x\).
With respect to the chosen invariant inner product, write \(\mathsf B(v,w)=\langle \mathcal A v,w\rangle\).
Since the self-adjoint operator \(\mathcal A\) commutes with \(G_x\), its negative and positive spectral subspaces are the invariant orthogonal summands \(V^-\) and \(V^+\).

By Hepworth's Morse lemma~\cite[Theorem~4.7 and Lemma~4.9]{hepworth2009morse}, there is a \(G_x\)-equivariant diffeomorphism \(\phi\) from a neighborhood of the origin onto a neighborhood of \(\widetilde x\) such that
\[
\widetilde f(\phi(u))=f(x)+\mathsf B(u,u).
\]
For each eigenvalue \(\lambda\) of \(\mathcal A\), let \(E_\lambda\) be the corresponding eigenspace and define \(C|_{E_\lambda}=|\lambda|^{-1/2}\Id_{E_\lambda}\).
The map \(C\) is \(G_x\)-equivariant and satisfies
\[
\mathsf B(Cu,Cu)
=
-\lVert u^-\rVert^2+\lVert u^+\rVert^2.
\]
Use the local coordinates induced by \(\phi\circ C\).
The Morse index is \(\dim V^-\).
\end{proof}

\begin{lemma}\label{lem:main-hessian-splitting-principle}
Let \(G<\Orth(n)\) be finite, and let \(\mathbb R^n=V^-\oplus V^+\) be an orthogonal \(G\)-invariant splitting.
Define \(Q_{V^-,V^+}\colon\mathbb R^n\to\mathbb R\) by
\[
Q_{V^-,V^+}(v^-+v^+)
=
-\lVert v^-\rVert^2+\lVert v^+\rVert^2,
\qquad v^-\in V^-,\quad v^+\in V^+.
\]
The function \(Q_{V^-,V^+}\) is \(G\)-invariant.
It has a nondegenerate critical point at the origin.
The negative and positive eigenspaces of its Hessian at the origin are \(V^-\) and \(V^+\), respectively.
Its Morse index at the origin is \(\dim V^-\).
\end{lemma}

\begin{proof}
The action of \(G\) preserves both summands and their norms, so \(Q_{V^-,V^+}\) is \(G\)-invariant.
The origin is a critical point.
In orthonormal bases of \(V^-\) and \(V^+\), the Hessian matrix is
\[
\operatorname{diag}(-2I_{\dim V^-},2I_{\dim V^+}).
\]
\end{proof}

For the splitting in Lemma~\ref{lem:orbifold-morse-normal-form}, let \(G<\Orth(n)\) be the image of the isotropy representation.
The associated quotient handle is \(D(V^-)\times_GD(V^+)\).
Its attaching and belt regions are
\[
S(V^-)\times_GD(V^+)
\qquad\text{and}\qquad
D(V^-)\times_GS(V^+),
\]
respectively.
Their common boundary is \(S(V^-)\times_GS(V^+)\).
When \(G\) is trivial and \(k=\dim V^-\), the quotient handle is diffeomorphic to the usual smooth handle \(D^k\times D^{n-k}\).
We use the term \textup{smooth handle} only when \(G\) is trivial.
In dimension four, let \(\rho\colon G\hookrightarrow\SO(4)\) be the inclusion.
The quotient handle is \(H_\rho(V^-)\) of Definition~\ref{def:ordered-handle}.

\subsection{Ordered handles and normalizer orbits}\label{subsec:main-quotient-handle}

Fix a finite group \(G\) and a faithful representation \(\rho\colon G\to\SO(4)\).
Let \(0\le k\le4\).
For \(V\in\operatorname{Gr}_k(\rho)\), the quotient handle \(H_\rho(V)\) alone need not determine the invariant splitting \(\mathbb R^4=V\oplus V^\perp\).
The attaching region, belt region, core, cocore, and origin in \(\mathbf H_\rho(V)\) record the two summands and their order.

Every \(T\in N^+_\rho\) induces an automorphism of \(G\) by conjugation.
The action of \(N^+_\rho\) on \(\operatorname{Gr}_k(\rho)\) is given by \(T\cdot V=T(V)\).

\begin{proposition}\label{prop:ordered-handle-orbits}
The assignment \(V\longmapsto\mathbf H_\rho(V)\) induces a bijection between the \(N^+_\rho\)-orbits in \(\operatorname{Gr}_k(\rho)\) and the isomorphism classes of ordered handles determined by \(\rho(G)\)-invariant \(k\)-dimensional subspaces.
The ordered dual of \(\mathbf H_\rho(V)\) is \(\mathbf H_\rho(V^\perp)\).
The ordered handle \(\mathbf H_\rho(V)\) is isomorphic to its ordered dual if and only if an element of \(N^+_\rho\) carries \(V\) to \(V^\perp\).
\end{proposition}

\begin{proof}
Let \(V,V'\in\operatorname{Gr}_k(\rho)\).
Suppose that \(S\in N^+_\rho\) satisfies \(S(V)=V'\).
Then \(S(V^\perp)=(V')^\perp\), and \(S\) induces an isomorphism
\[
\mathbf H_\rho(V)\longrightarrow\mathbf H_\rho(V').
\]

Conversely, let \(F\colon \mathbf H_\rho(V)\to\mathbf H_\rho(V')\) be an isomorphism.
Choose a local lift \(\widetilde F\) near the origin with \(\widetilde F(0)=0\).
Since \(F\) is an orbifold diffeomorphism, \(d\widetilde F_0\) is invertible.
The homomorphism of local groups associated with the lift \(\widetilde F\) is an automorphism.
Since \(\rho\) is faithful, this automorphism has the form
\[
\rho(g)\longmapsto\rho(\alpha(g))
\]
for a unique \(\alpha\in\operatorname{Aut}(G)\).
Set \(L=d\widetilde F_0\).
Equivariance gives
\[
L\rho(g)=\rho(\alpha(g))L
\]
for every \(g\in G\).

Under the quotient projections, the preimages of the core and cocore in the source are \(D(V)\times\{0\}\) and \(\{0\}\times D(V^\perp)\), respectively.
The preimages of the core and cocore in the target are \(D(V')\times\{0\}\) and \(\{0\}\times D((V')^\perp)\), respectively.
Since \(F\) preserves the core and cocore,
\[
L(V)=V',
\qquad
L(V^\perp)=(V')^\perp.
\]

Write the polar decomposition as \(L=TP\), where \(P=(L^*L)^{1/2}\) is positive and \(T\) is orthogonal.
The equivariance relation implies that \(L^*L\) commutes with \(\rho(G)\).
Its positive square root \(P\) also commutes with \(\rho(G)\).
For \(v\in V\) and \(w\in V^\perp\), we have \(Lv\in V'\) and \(Lw\in(V')^\perp\).
Therefore
\[
\langle L^*Lv,w\rangle=\langle Lv,Lw\rangle=0.
\]
Thus \(L^*L\) preserves \(V\) and \(V^\perp\).
Its positive square root \(P\) preserves both subspaces.
Since \(\widetilde F\) preserves orientation, \(T=LP^{-1}\) lies in \(\SO(4)\).
Moreover,
\[
T(V)=V',
\qquad
T\rho(g)T^{-1}=\rho(\alpha(g))
\]
for every \(g\in G\).
Thus \(T\in N^+_\rho\).

The ordered dual of \(\mathbf H_\rho(V)\) is \(\mathbf H_\rho(V^\perp)\).
If \(k\ne2\), then \(V\) and \(V^\perp\) have different dimensions, so no element of \(N^+_\rho\) carries \(V\) to \(V^\perp\).
For \(k=2\), take \(V'=V^\perp\).
\end{proof}

\begin{proof}[\textup{Proof of Theorem A}]
By Proposition~\ref{prop:ordered-handle-orbits}, the \(N^+_\rho\)-orbits in \(\operatorname{Gr}_k(\rho)\) are in bijection with the isomorphism classes of ordered handles.
By Lemma~\ref{lem:main-hessian-splitting-principle}, the invariant quadratic function \(Q_{V,V^\perp}\) realizes the isomorphism class of ordered handles represented by \(V\).
By Definition~\ref{def:ordered-handle-isomorphism}, the ordered dual of \(\mathbf H_\rho(V)\) is \(\mathbf H_\rho(V^\perp)\).
Every \(T\in N^+_\rho\) is orthogonal.
For every \(\rho(G)\)-invariant subspace \(V\), we therefore have \(T(V^\perp)=T(V)^\perp\).
Thus orthogonal complementation is well defined on the \(N^+_\rho\)-orbits of \(\rho(G)\)-invariant subspaces of \(\mathbb R^4\).
\end{proof}

\section{Invariant splittings and Morse indices}\label{sec:morse-index-classification}

In this section, we determine the possible dimensions of \(G\)-invariant subspaces of \(\mathbb R^4\) for each finite subgroup \(G<\SO(4)\).
By Lemmas~\ref{lem:orbifold-morse-normal-form} and~\ref{lem:main-hessian-splitting-principle}, these dimensions are precisely the possible Morse indices.
Section~\ref{subsec:so4-finite-subgroups} introduces the quaternionic representatives for the families, fixes the parameter conventions, and defines factor interchange.
Section~\ref{subsec:q8-criterion} proves criteria for irreducibility and determines the invariant subspaces of the linear transformations used to calculate the possible indices.
Section~\ref{subsec:parameter-cases} determines the parameter conditions for the existence of invariant lines in \textup{Families} \(1\) and \(1'\).
It also determines the conditions for the existence of invariant \(2\)-dimensional subspaces other than \(\operatorname{span}\{1,\mathbf{i}\}\) and \(\operatorname{span}\{\mathbf{j},\mathbf{k}\}\).
Section~\ref{subsec:morse-index-theorem} applies these criteria to classify the possible Morse indices and prove Theorem~\ref{thm:4d-morse-index-classification}.

\subsection{Finite subgroups of \texorpdfstring{\(\SO(4)\)}{SO(4)} and the completed Du Val list}\label{subsec:so4-finite-subgroups}

Threlfall and Seifert classified the finite subgroups of \(\SO(4)\) up to conjugacy~\cite{seifertthrelfall1931,seifertthrelfall1933}.
We use Du Val's quaternionic notation~\cite{val1964homographies} and the completion of Du Val's list by Mecchia and Seppi.
Table~\ref{tab:mecchia-seppi-families} lists the completed Du Val families.
The subscripts \(1\), \(1'\), and \(26''\) are family labels in the completed Du Val list.
For \(q\ge1\), let \(\Cyc{q}\) denote the cyclic group of order \(q\).
Write \(\mathbf T,\mathbf O,\mathbf I<\SO(3)\) for the tetrahedral, octahedral, and icosahedral rotation groups, respectively.
The binary lift of a finite rotation group is its full inverse image under the standard double covering from the group of unit quaternions to \(\SO(3)\).
This covering sends a unit quaternion \(q\) to the rotation \(v\mapsto qvq^{-1}\) on the imaginary quaternions.
The symbols \(T^*,O^*,I^*\) denote the binary tetrahedral, binary octahedral, and binary icosahedral groups, respectively.
For \(q\ge1\), let \(\Dih{2q}<\SO(3)\) denote the rotational dihedral group of order \(2q\).
For \(q\ge2\), let \(\BDih{4q}\) denote its binary lift, which has order \(4q\).
We use the convention \(\Dih{2}=\Cyc{2}\) when \(q=1\).
We restrict the notation \(\BDih{4q}\) to \(q\ge2\), as in Mecchia and Seppi~\cite[p.~814]{mecchia2015fibered}.
The smallest binary dihedral group in this notation is \(\BDih{8}\).
The subscript \(k\) in \(\BDih{k}\) is the order of the group.

Let
\[
\mathbb H=\{a+b\mathbf{i}+c\mathbf{j}+d\mathbf{k}\mid a,b,c,d\in\mathbb R\}
=\{z_1+z_2\mathbf{j}\mid z_1,z_2\in\mathbb C\}
\]
be the quaternion algebra.
Let \(S^3=\{q\in\mathbb H\mid \lVert q\rVert=1\}\) be the group of unit quaternions.
We realize \(\Cyc{q}\) as the subgroup \(\langle e^{2\pi \mathbf{i}/q}\rangle<S^3\).
We regard \(T^*\), \(O^*\), \(I^*\), and \(\BDih{4q}\) for \(q\ge2\) as finite subgroups of \(S^3\).

Define
\[
\Phi\colon S^3\times S^3\longrightarrow\SO(4),
\qquad
\Phi(p,q)(v)=pvq^{-1},
\quad v\in\mathbb H.
\]
For \(i=1,2\), let \(\pi_i\colon S^3\times S^3\to S^3\) be the projection to the \(i\)-th factor.
Finite subgroups of \(\SO(4)\) admit the following quaternionic description~\cite[Section~2]{mecchia2015fibered}.
Let \(\widetilde G<S^3\times S^3\) be finite.
Set \(L=\pi_1(\widetilde G)\) and \(R=\pi_2(\widetilde G)\).
Set
\[
L_K=\pi_1\bigl(\widetilde G\cap(S^3\times\{1\})\bigr)
\]
and
\[
R_K=\pi_2\bigl(\widetilde G\cap(\{1\}\times S^3)\bigr).
\]

The following proposition follows from the description of finite subgroups of \(\SO(4)\) in~\cite[Section~2, p.~813]{mecchia2015fibered}.

\begin{proposition}\label{prop:quaternionic-cover-fiber-product}
The following hold.
\begin{enumerate}[label=\textup{(\roman*)}]
\item
The map \(\Phi\) is a two-fold covering homomorphism with
\(\ker\Phi=\{(1,1),(-1,-1)\}\).
The finite subgroups of \(\SO(4)\) correspond to the finite subgroups
of \(S^3\times S^3\) containing \(\ker\Phi\).

\item
The subgroups \(L_K\) and \(R_K\) are normal in \(L\) and \(R\),
respectively.
There is a unique isomorphism
\(\varphi_{\widetilde G}\colon L/L_K\to R/R_K\)
such that \(\varphi_{\widetilde G}(lL_K)=rR_K\)
for every \((l,r)\in\widetilde G\).
Moreover,
\[
\widetilde G
=
\{(l,r)\in L\times R
\mid\varphi_{\widetilde G}(lL_K)=rR_K\}.
\]
\end{enumerate}
\end{proposition}

Let \(L_K\triangleleft L\) and \(R_K\triangleleft R\).
Let \(\varphi\colon L/L_K\to R/R_K\) be an isomorphism.
We write
\[
(L/L_K,R/R_K)_\varphi
=
\{(l,r)\in L\times R\mid
\varphi(lL_K)=rR_K\}.
\]
If the conjugacy class of \((L/L_K,R/R_K)_\varphi\) is determined by \(L\), \(L_K\), \(R\), and \(R_K\), we omit the subscript \(\varphi\), as in the completed Du Val list.
Table~\ref{tab:mecchia-seppi-families} lists the completed Du Val families of Mecchia and Seppi~\cite[Table~1]{mecchia2015fibered}.
Its entries describe subgroups of \(S^3\times S^3\), and the order column records the orders of their images under \(\Phi\).
The integers \(m,n,r\) are positive whenever they occur.
The parameter \(s\) is an integer.
Table~\ref{tab:mecchia-seppi-seifert-primed} lists the Seifert invariants for \textup{Families} \(1'\) and \(11'\).
Table~\ref{tab:mecchia-seppi-seifert-unprimed} lists those for \textup{Families} \(1\) and \(11\), as computed by Mecchia and Seppi~\cite[Tables~2 and~3]{mecchia2015fibered}.
The notation in each table is local to that table.

For a \(2\)-manifold \(X\) without boundary, \(X(n_1,\ldots,n_k)\) denotes the \(2\)-orbifold with underlying space \(X\) and cone points of indices \(n_1,\ldots,n_k\)~\cite[p.~817]{mecchia2015fibered}.
If \(X\) has nonempty connected boundary, then \(X(n_1,\ldots,n_k;m_1,\ldots,m_h)\) denotes the \(2\)-orbifold with underlying space \(X\), cone points of indices \(n_1,\ldots,n_k\), and corner reflectors of indices \(m_1,\ldots,m_h\).

\begin{table}[p]
\centering
\fontsize{9}{11.4}\selectfont
\setlength{\tabcolsep}{6pt}
\setlength{\arrayrulewidth}{0.7pt}
\renewcommand{\arraystretch}{1}
\begin{tabular}{|l|c|c|c|}
\hline
 & family of groups & order of \(\Phi(G)\) & \\
\hline
\(1\). & \((C_{2mr}/C_{2m},C_{2nr}/C_{2n})_s\) & \(2mnr\) & \(\gcd(s,r)=1\)\\
\(1'\). & \((C_{mr}/C_{m},C_{nr}/C_{n})_s\) & \((mnr)/2\) & \(\gcd(s,r)=1\,\gcd(2,n)=1\)\\
 & & & \(\gcd(2,m)=1\,\gcd(2,r)=2\)\\
\(2\). & \((C_{2m}/C_{2m},D_{4n}^{*}/D_{4n}^{*})\) & \(4mn\) & \\
\(3\). & \((C_{4m}/C_{2m},D_{4n}^{*}/C_{2n})\) & \(4mn\) & \\
\(4\). & \((C_{4m}/C_{2m},D_{8n}^{*}/D_{4n}^{*})\) & \(8mn\) & \\
\(5\). & \((C_{2m}/C_{2m},T^*/T^*)\) & \(24m\) & \\
\(6\). & \((C_{6m}/C_{2m},T^*/D_{8}^{*})\) & \(24m\) & \\
\(7\). & \((C_{2m}/C_{2m},O^*/O^*)\) & \(48m\) & \\
\(8\). & \((C_{4m}/C_{2m},O^*/T^*)\) & \(48m\) & \\
\(9\). & \((C_{2m}/C_{2m},I^*/I^*)\) & \(120m\) & \\
\(10\). & \((D_{4m}^{*}/D_{4m}^{*},D_{4n}^{*}/D_{4n}^{*})\) & \(8mn\) & \\
\(11\). & \((D_{4mr}^{*}/C_{2m},D_{4nr}^{*}/C_{2n})_s\) & \(4mnr\) & \(\gcd(s,r)=1\)\\
\(11'\). & \((D_{2mr}^{*}/C_{m},D_{2nr}^{*}/C_{n})_s\) & \(mnr\) & \(\gcd(s,r)=1\,\gcd(2,n)=1\)\\
 & & & \(\gcd(2,m)=1\,\gcd(2,r)=2\)\\
\(12\). & \((D_{8m}^{*}/D_{4m}^{*},D_{8n}^{*}/D_{4n}^{*})\) & \(16mn\) & \\
\(13\). & \((D_{8m}^{*}/D_{4m}^{*},D_{4n}^{*}/C_{2n})\) & \(8mn\) & \\
\(14\). & \((D_{4m}^{*}/D_{4m}^{*},T^*/T^*)\) & \(48m\) & \\
\(15\). & \((D_{4m}^{*}/D_{4m}^{*},O^*/O^*)\) & \(96m\) & \\
\(16\). & \((D_{4m}^{*}/C_{2m},O^*/T^*)\) & \(48m\) & \\
\(17\). & \((D_{8m}^{*}/D_{4m}^{*},O^*/T^*)\) & \(96m\) & \\
\(18\). & \((D_{12m}^{*}/C_{2m},O^*/D_{8}^{*})\) & \(48m\) & \\
\(19\). & \((D_{4m}^{*}/D_{4m}^{*},I^*/I^*)\) & \(240m\) & \\
\(20\). & \((T^*/T^*,T^*/T^*)\) & \(288\) & \\
\(21\). & \((T^*/C_{2},T^*/C_{2})\) & \(24\) & \\
\(21'\). & \((T^*/C_{1},T^*/C_{1})\) & \(12\) & \\
\(22\). & \((T^*/D_{8}^{*},T^*/D_{8}^{*})\) & \(96\) & \\
\(23\). & \((T^*/T^*,O^*/O^*)\) & \(576\) & \\
\(24\). & \((T^*/T^*,I^*/I^*)\) & \(1440\) & \\
\(25\). & \((O^*/O^*,O^*/O^*)\) & \(1152\) & \\
\(26\). & \((O^*/C_{2},O^*/C_{2})\) & \(48\) & \\
\(26'\). & \((O^*/C_{1},O^*/C_{1})_{Id}\) & \(24\) & \\
\(26''\). & \((O^*/C_{1},O^*/C_{1})_f\) & \(24\) & \\
\(27\). & \((O^*/D_{8}^{*},O^*/D_{8}^{*})\) & \(192\) & \\
\(28\). & \((O^*/T^*,O^*/T^*)\) & \(576\) & \\
\(29\). & \((O^*/O^*,I^*/I^*)\) & \(2880\) & \\
\(30\). & \((I^*/I^*,I^*/I^*)\) & \(7200\) & \\
\(31\). & \((I^*/C_{2},I^*/C_{2})_{Id}\) & \(120\) & \\
\(31'\). & \((I^*/C_{1},I^*/C_{1})_{Id}\) & \(60\) & \\
\(32\). & \((I^*/C_{2},I^*/C_{2})_f\) & \(120\) & \\
\(32'\). & \((I^*/C_{1},I^*/C_{1})_f\) & \(60\) & \\
\(33\). & \((D_{8m}^{*}/C_{2m},D_{8n}^{*}/C_{2n})_f\) & \(8mn\) & \(m\ne1\quad n\ne1\).\\
\(33'\). & \((D_{8m}^{*}/C_{m},D_{8n}^{*}/C_{n})_f\) & \(4mn\) & \(\gcd(2,n)=1\,\gcd(2,m)=1\)\\
 & & & \(m\ne1\text{ and }n\ne1\).\\
\(34\). & \((C_{4m}/C_{m},D_{4n}^{*}/C_{n})\) & \(2mn\) & \(\gcd(2,n)=1\,\gcd(2,m)=1\)\\
\hline
\end{tabular}
\vspace{12pt}
\caption{Finite subgroups of \(\SO(4)\), from Mecchia and Seppi~\cite[Table~1]{mecchia2015fibered}.}\label{tab:mecchia-seppi-families}\label{tab:mecchia-seppi-families-1-19}\label{tab:mecchia-seppi-families-20-34}
\end{table}

\begin{table}[p]
\centering
\normalsize
\setlength{\tabcolsep}{6pt}
\renewcommand{\arraystretch}{1}
\begin{tabular}{|p{\dimexpr\linewidth-2\tabcolsep-2\arrayrulewidth\relax}|}
\hline
\begin{minipage}{\linewidth}
\centering
\vspace{9pt}
\(G=\Phi((C_{mr}/C_m,C_{nr}/C_n)_s)\) (\textup{Family} \(1'\))\par
\vspace{3pt}
and\par
\vspace{3pt}
\(G=\Phi((D_{2mr}^{*}/C_m,D_{2nr}^{*}/C_n)_s)\) (\textup{Family} \(11'\))\par
\vspace{10pt}
\end{minipage}\\
\hline
\begin{minipage}{\linewidth}
\vspace{9pt}
\begin{tabular}{@{}l@{\quad}l@{}}
We define: & \(h=\gcd(m,n),\ m'=\frac{m}{h},\ n'=\frac{n}{h},\)\\[3pt]
 & \(a=\gcd(n'-sm',n'+sm',m'n'r),\)\\[3pt]
 & \(b_1=\gcd\bigl(\frac{n'-sm'}{a},\frac{m'n'r}{a}\bigr),\)\\[3pt]
 & \(b_2=\gcd\bigl(\frac{n'+sm'}{a},\frac{m'n'r}{a}\bigr),\)\\[3pt]
 & \(\nu\) minimal positive integer s.t.\ \(\gcd\bigl(\frac{2n'}{a\nu},\frac{a}{2}\bigr)=1,\)\\[3pt]
 & \(d=\frac{\nu^2a(n'+sm')+2n'm'r}{2a\nu b_2},\)\\[3pt]
 & \(g=\frac{\nu^2a(n'-sm')-2n'm'r}{2a\nu b_1},\)\\[3pt]
 & \(f=\frac{m'n'r}{2b_1b_2},\)\\[3pt]
 & \(\bar g\) such that \(g\bar g\equiv1\bmod f,\)\\[3pt]
 & \(\bar c\) such that \(\bigl(\nu s+r\frac{2n'}{a\nu}\bigr)\bar c\equiv1\bmod n'r.\)
\end{tabular}\par
\vspace{10pt}
\end{minipage}\\
\hline
\begin{minipage}{\linewidth}
\setlength{\parindent}{0pt}
\setlength{\parskip}{6pt}
\vspace{6pt}
The orbifold \(S^3/\Phi((C_{mr}/C_m,C_{nr}/C_n)_s)\) fibers over \(S^2(nr/2,nr/2)\) with local invariants \((d\bar c b_2h)/(nr/2)\) and \(-(g\bar c b_1h)/(nr/2)\) and Euler number \(-2m/(nr)\).\par
The underlying topological space of \(S^3/\Phi((C_{mr}/C_m,C_{nr}/C_n)_s)\) is the lens space \(L(f,d\bar g)\).\par
The singular set of \(S^3/\Phi((C_{mr}/C_m,C_{nr}/C_n)_s)\) is a link with at most two components of singular index \(b_2h\) and \(b_1h\) (if the singular index is \(1\) the corresponding component consists of non-singular points).\par
\vspace{5pt}
\end{minipage}\\
\hline
\begin{minipage}{\linewidth}
\setlength{\parindent}{0pt}
\setlength{\parskip}{9pt}
\vspace{6pt}
The orbifold \(S^3/\Phi((D_{2mr}^{*}/C_m,D_{2nr}^{*}/C_n)_s)\) fibers over \(D^2(;nr/2,nr/2)\) with local invariants \((d\bar c b_2h)/(nr/2)\) and \(-(g\bar c b_1h)/(nr/2)\) and Euler number \(-m/(nr)\).\par
The underlying topological space of \(S^3/\Phi((D_{2mr}^{*}/C_m,D_{2nr}^{*}/C_n)_s)\) is the 3-sphere.\par
\vspace{6pt}
\end{minipage}\\
\hline
\end{tabular}
\vspace{8pt}
\caption{\textup{Families} \(1'\) and \(11'\), from Mecchia and Seppi~\cite[Table~2]{mecchia2015fibered}. In the definition of \(a\), we replace \(m'+sn'\) in the original table with \(n'+sm'\), as in~\cite[p.~826]{mecchia2015fibered}.}\label{tab:mecchia-seppi-seifert-primed}
\end{table}

\begin{table}[p]
\centering
\normalsize
\setlength{\tabcolsep}{6pt}
\renewcommand{\arraystretch}{1}
\begin{tabular}{|p{\dimexpr .60\linewidth-2\tabcolsep-1.5\arrayrulewidth\relax}|p{\dimexpr .40\linewidth-2\tabcolsep-1.5\arrayrulewidth\relax}|}
\hline
\multicolumn{2}{|p{\dimexpr\linewidth-2\tabcolsep-2\arrayrulewidth\relax}|}{%
\begin{minipage}{\linewidth}
\centering
\vspace{9pt}
\(G=\Phi((C_{2mr}/C_{2m},C_{2nr}/C_{2n})_s)\) (\textup{Family} \(1\))\par
\vspace{3pt}
and\par
\vspace{3pt}
\(G=\Phi((D^*_{4mr}/C_{2m},D^*_{4nr}/C_{2n})_s)\) (\textup{Family} \(11\))\par
\vspace{10pt}
\end{minipage}
}\tabularnewline
\hline
\multicolumn{2}{|p{\dimexpr\linewidth-2\tabcolsep-2\arrayrulewidth\relax}|}{%
\begin{minipage}{\linewidth}
\setlength{\parindent}{0pt}
\vspace{9pt}
We define: \(h=\gcd(m,n),\ m'=\frac{m}{h},\ n'=\frac{n}{h},\)\par
\vspace{3pt}
\(a=\gcd(n'-sm',n'+sm',2m'n'r),\)\par
\vspace{3pt}
\(b_1=\gcd\bigl(\frac{n'-sm'}{a},\frac{2m'n'r}{a}\bigr),\ b_2=\gcd\bigl(\frac{n'+sm'}{a},\frac{2m'n'r}{a}\bigr).\)\par
\vspace{9pt}
Remark: without loss of generality, we assume \(s\) odd.\par
\vspace{10pt}
\end{minipage}
}\tabularnewline
\hline
\begin{minipage}[t]{\linewidth}
\setlength{\parindent}{0pt}
\vspace{9pt}
If \(n'm'\) is even, we define:\par
\vspace{5pt}
\(\nu\) minimal positive integer s.t.\ \(\gcd\bigl(\frac{n'}{a\nu},a\bigr)=1,\)\par
\vspace{5pt}
\(f_1=f_2=1.\)\par
\vspace{15pt}
If \(n'm'\) is odd, we define:\par
\vspace{5pt}
\(\nu\) minimal positive integer s.t.\ \(\gcd\bigl(\frac{2n'}{a\nu},\frac a2\bigr)=1,\)\par
\vspace{5pt}
\(f_i=\begin{cases}
2&\text{if }r/b_i\text{ is even},\\
1&\text{if }r/b_i\text{ is odd}.
\end{cases}\)\par
\vspace{10pt}
\end{minipage}
&
\begin{minipage}[t]{\linewidth}
\setlength{\parindent}{0pt}
\vspace{9pt}
In both cases we define:\par
\vspace{5pt}
\(d=\frac{\nu^2a(n'+sm')+2n'm'r}{f_2a\nu b_2},\)\par
\vspace{5pt}
\(g=\frac{\nu^2a(n'-sm')-2n'm'r}{f_1a\nu b_1},\)\par
\vspace{5pt}
\(f=\frac{2m'n'r}{f_1f_2b_1b_2},\)\par
\vspace{5pt}
\(\bar g\) s.t.\ \(g\bar g\equiv1\bmod f,\)\par
\vspace{5pt}
\(\bar c\) s.t.\par
\vspace{8pt}
\(\bigl(\nu s+r\frac{2n'}{a\nu}\bigr)\bar c\equiv1\bmod n'r.\)\par
\vspace{10pt}
\end{minipage}
\tabularnewline
\hline
\multicolumn{2}{|p{\dimexpr\linewidth-2\tabcolsep-2\arrayrulewidth\relax}|}{%
\begin{minipage}{\linewidth}
\setlength{\parindent}{0pt}
\setlength{\parskip}{6pt}
\vspace{6pt}
The orbifold \(S^3/\Phi((C_{2mr}/C_{2m},C_{2nr}/C_{2n})_s)\) fibers over \(S^2(nr,nr)\) with local invariants \(d\bar c f_2b_2h/(nr)\) and \(-g\bar c f_1b_1h/(nr)\) and Euler number \(-2m/(nr)\).\par
The underlying topological space of \(S^3/\Phi((C_{2mr}/C_{2m},C_{2nr}/C_{2n})_s)\) is the lens space \(L(f,d\bar g)\).\par
The singular set of \(S^3/\Phi((C_{2mr}/C_{2m},C_{2nr}/C_{2n})_s)\) is a link with at most two components of singular index \(f_2b_2h\) and \(f_1b_1h\) (if the singular index is \(1\) the corresponding component consists of non-singular points).\par
\vspace{5pt}
\end{minipage}
}\tabularnewline
\hline
\multicolumn{2}{|p{\dimexpr\linewidth-2\tabcolsep-2\arrayrulewidth\relax}|}{%
\begin{minipage}{\linewidth}
\setlength{\parindent}{0pt}
\setlength{\parskip}{9pt}
\vspace{6pt}
The orbifold \(S^3/\Phi((D^*_{4mr}/C_{2m},D^*_{4nr}/C_{2n})_s)\) fibers over \(D^2(;nr,nr)\) with local invariants \(d\bar c f_2b_2h/(nr)\) and \(-g\bar c f_1b_1h/(nr)\) and Euler number \(-m/(nr)\).\par
The underlying topological space of \(S^3/\Phi((D^*_{4mr}/C_{2m},D^*_{4nr}/C_{2n})_s)\) is the 3-sphere.\par
\vspace{6pt}
\end{minipage}
}\tabularnewline
\hline
\end{tabular}
\vspace{8pt}
\caption{\textup{Families} \(1\) and \(11\), from Mecchia and Seppi~\cite[Table~3]{mecchia2015fibered}. In the definition of \(a\), we replace \(m'+sn'\) in the original table with \(n'+sm'\), as in~\cite[p.~830]{mecchia2015fibered}.}\label{tab:mecchia-seppi-seifert-unprimed}
\end{table}

Let \(\epsilon_{26}\colon O^*\longrightarrow\{\pm1\}\) be the homomorphism with kernel \(T^*\).
Define \(f_0\colon O^*\longrightarrow O^*\) by
\[
f_0(q)=\epsilon_{26}(q)q.
\]
The automorphism \(f_0\) represents the nontrivial outer automorphism class of \(O^*\)~\cite[p.~57]{val1964homographies}.

The subscript \(s\) specifies the quotient isomorphism \(\varphi_s\) in Table~\ref{tab:mecchia-seppi-quotient-isomorphisms}.
For each family, the images in Table~\ref{tab:mecchia-seppi-quotient-isomorphisms} determine the quotient isomorphism.
The subscript \(f\) denotes the quotient isomorphism for the family in which it occurs and need not have the same meaning in different families.
For \textup{Family} \(26''\), the quotient isomorphism is the automorphism \(f_0\).
For \textup{Families} \(32\) and \(32'\), the subscript \(f\) represents the nontrivial outer automorphism class of \(I^*/\Cyc2\) and \(I^*\), respectively.
Section~\ref{subsec:appendix-polyhedral} fixes the chosen representatives for these two families.
The tables specify conjugacy classes.
We fix representatives in quaternionic coordinates for \textup{Families} \(26''\), \(32\), and \(32'\) below and in Section~\ref{subsec:appendix-polyhedral}.

\begin{table}[!htbp]
\centering
\small
\renewcommand{\arraystretch}{1.6}
\begin{tabular}{@{}ll@{}}
\hline
\textup{Family} & Quotient isomorphism\\
\hline
\(1\), \(11\) & \(\varphi_s(e^{\pi\mathbf{i}/(mr)}\Cyc{2m})=e^{s\pi\mathbf{i}/(nr)}\Cyc{2n}\)\\
\(11\) & \(\varphi_s(\mathbf{j}\Cyc{2m})=\mathbf{j}\Cyc{2n}\)\\
\(1'\), \(11'\) & \(\varphi_s(e^{2\pi\mathbf{i}/(mr)}\Cyc m)=e^{2s\pi\mathbf{i}/(nr)}\Cyc n\)\\
\(11'\) & \(\varphi_s(\mathbf{j}\Cyc m)=\mathbf{j}\Cyc n\)\\
\(33\) & \(f(e^{\pi\mathbf{i}/(2m)}\Cyc{2m})=\mathbf{j}\Cyc{2n}\)\\
 & \(f(\mathbf{j}\Cyc{2m})=e^{\pi\mathbf{i}/(2n)}\Cyc{2n}\)\\
\(33'\) & \(f(e^{\pi\mathbf{i}/(2m)}\Cyc m)=\mathbf{j}\Cyc n\)\\
 & \(f(\mathbf{j}\Cyc m)=e^{\pi\mathbf{i}/(2n)}\Cyc n\)\\
\hline
\end{tabular}
\caption{Quotient isomorphisms for the completed Du Val list~\cite[pp.~814--816]{mecchia2015fibered}. The row for \textup{Family} \(33'\) includes the correction in Remark~\ref{rem:family-33p-angle-correction}.}\label{tab:mecchia-seppi-quotient-isomorphisms}
\end{table}

For each family label \(\nu\ne26''\), let \(\widetilde G_\nu<S^3\times S^3\) be the representative specified in the completed Du Val list.
Set \(G_\nu=\Phi(\widetilde G_\nu)\).

For \textup{Family} \(26''\), choose
\[
\widetilde G_{26''}
=
\{(q,\mathbf{i}^{-1}f_0(q)\mathbf{i})\mid q\in O^*\},
\qquad
G_{26''}=\Phi(\widetilde G_{26''}).
\]

When a family is fixed, we write \(G\) for its chosen representative.
Let \(\rho\colon G\hookrightarrow\SO(4)\) denote the inclusion.

\begin{definition}\label{def:representative-family-parameters}
For the chosen representatives, the parameters satisfy the following conditions.
Every integer parameter other than \(s\) is positive whenever it occurs.
We allow \(s\in\mathbb Z\).
For \textup{Families} \(1\) and \(11\), assume that \(\gcd(s,r)=1\).
For \textup{Families} \(1'\) and \(11'\), assume that \(m\) and \(n\) are odd.
Write \(r=2r'\) for a positive integer \(r'\), and assume that \(\gcd(s,2r')=1\).
The additional parameter ranges are
\[
\begin{array}{c|c}
\textnormal{families}&\textnormal{additional ranges}\\
\hline
2,3,4&m\ge1,\ n\ge2\\
10,12,13&m,n\ge2\\
11&mr\ge2,\ nr\ge2\\
11'&mr'\ge2,\ nr'\ge2\\
14,15,16,17,19&m\ge2\\
33&m,n\ge2\\
33'&m,n\ge3,\quad m,n\text{ odd}\\
34&m\ge1,\quad n\ge3,\quad m,n\text{ odd}.
\end{array}
\]
For \textup{Families} \(5\)--\(9\) and \(18\), we assume that \(m\ge1\).
\end{definition}

Define quaternionic conjugation \(\mu\colon\mathbb H\longrightarrow\mathbb H\) by \(\mu(v)=\overline v\).
For every \(q\in S^3\), we have \(\mu(q)=q^{-1}\).
With respect to the basis \((1,\mathbf{i},\mathbf{j},\mathbf{k})\), the matrix of \(\mu\) is
\[
\operatorname{diag}(1,-1,-1,-1).
\]
This matrix has determinant \(-1\), so \(\mu\in\Orth(4)\setminus\SO(4)\).

For \(p,q\in S^3\),
\[
\mu\Phi(p,q)\mu^{-1}=\Phi(q,p).
\]

\begin{definition}\label{def:factor-interchange}
Let \(\widetilde G<S^3\times S^3\) be a finite subgroup containing \(\ker\Phi\).
Set \(G=\Phi(\widetilde G)\).
Define
\[
G^\mu=\mu G\mu^{-1}
\]
and
\[
\widetilde G^\mu=\{(r,l)\mid(l,r)\in\widetilde G\}.
\]
We have \(G^\mu=\Phi(\widetilde G^\mu)\).
The assignment \(\widetilde G\mapsto\widetilde G^\mu\), or equivalently \(G\mapsto G^\mu\), is called \emph{factor interchange}.
\end{definition}

The groups \(G\) and \(G^\mu\) are conjugate in \(\Orth(4)\) by the orientation-reversing isometry \(\mu\), but need not be conjugate in \(\SO(4)\).
Factor interchange exchanges \(L\) with \(R\) and \(L_K\) with \(R_K\).
It replaces \(\varphi_{\widetilde G}\) by \(\varphi_{\widetilde G}^{-1}\).
A subspace \(V\subset\mathbb H\) is \(G\)-invariant if and only if \(\mu V\) is \(G^\mu\)-invariant.
The map
\[
V^-\oplus V^+
\longmapsto
\mu V^-\oplus\mu V^+
\]
is a dimension-preserving bijection from the \(G\)-invariant orthogonal splittings to the \(G^\mu\)-invariant orthogonal splittings.
Conjugation by \(\mu\) defines an isomorphism
\[
N_{\SO(4)}(G)
\longrightarrow
N_{\SO(4)}(G^\mu),
\qquad
T\longmapsto\mu T\mu^{-1}.
\]
The bijection of invariant orthogonal splittings is equivariant for this isomorphism.
By Proposition~\ref{prop:ordered-handle-orbits}, \(G\) and \(G^\mu\) have the same number of isomorphism classes of ordered handles in each index.

Definition~\ref{def:representative-family-parameters} fixes the parameter ranges for the chosen representatives.
Each excluded choice of parameters determines a group that is \(\SO(4)\)-conjugate to one of these representatives or to its factor interchange~\cite[Remark~4.2 and Section~5.2]{mecchiaseppi2020diffeomorphism}.
Every \(\SO(4)\)-conjugacy class of finite subgroups of \(\SO(4)\) contains one of these representatives or its factor interchange~\cite[p.~816]{mecchia2015fibered}.

\begin{proposition}\label{prop:generator}
Let \(\widetilde G<S^3\times S^3\) be finite.
Set
\[
L=\pi_1(\widetilde G),
\qquad
R=\pi_2(\widetilde G).
\]
Set
\[
L_K=\pi_1\bigl(\widetilde G\cap(S^3\times\{1\})\bigr),
\qquad
R_K=\pi_2\bigl(\widetilde G\cap(\{1\}\times S^3)\bigr).
\]
Let \(a_1,\ldots,a_u\) generate \(L_K\), and let \(b_1,\ldots,b_v\) generate \(R_K\).
Choose generators \(\xi_1,\ldots,\xi_w\) of \(L/L_K\).
For \(1\le h\le w\), choose \((l_h,r_h)\in\widetilde G\) such that \(l_hL_K=\xi_h\).
Then \(\widetilde G\) is generated by
\[
(a_1,1),\ldots,(a_u,1),
\qquad
(1,b_1),\ldots,(1,b_v),
\qquad
(l_1,r_1),\ldots,(l_w,r_w).
\]
\end{proposition}

\begin{proof}
Let \(H\) be the subgroup generated by the elements in the proposition.
We have \(H\subset\widetilde G\), and the map \(H\to L/L_K\) is surjective.
For \((l,r)\in\widetilde G\), choose \((l',r')\in H\) such that \(l'L_K=lL_K\).
The isomorphism between the quotient groups gives \(r'R_K=rR_K\).
Thus
\[
(l',r')^{-1}(l,r)\in L_K\times R_K.
\]
Since \(L_K\times R_K\subset H\), we have \((l',r')^{-1}(l,r)\in H\).
It follows that \((l,r)\in H\).
\end{proof}

\clearpage
\subsection{Criteria for invariant subspaces}\label{subsec:q8-criterion}

Let \(Q_8=\{\pm1,\pm \mathbf{i},\pm \mathbf{j},\pm \mathbf{k}\}\subset S^3\subset\mathbb H\).

\begin{proposition}\label{prop:irred-local-gps}
Let \(\widetilde G<S^3\times S^3\) be a finite subgroup containing \(\ker\Phi\).
Write \(G=\Phi(\widetilde G)<\SO(4)\), and set
\[
L_K=\pi_1\bigl(\widetilde G\cap(S^3\times\{1\})\bigr)
\]
and
\[
R_K=\pi_2\bigl(\widetilde G\cap(\{1\}\times S^3)\bigr).
\]
Assume that either \(L_K\) or \(R_K\) generates \(\mathbb H\) as a real algebra.
Then the representation of \(G\) on \(\mathbb H\) is irreducible.
Every \(G\)-invariant smooth function on \(\mathbb R^4\) with a nondegenerate critical point at the origin has Morse index \(0\) or \(4\).
\end{proposition}

\begin{proof}
Let \(U\subset\mathbb H\) be a nonzero \(G\)-invariant subspace.
Choose \(0\ne u\in U\).

If \(L_K\) generates \(\mathbb H\), then \(U\) is invariant under left multiplication by every quaternion.
For every \(y\in\mathbb H\), we have \(y=(yu^{-1})u\in U\).
Hence \(U=\mathbb H\).

If \(R_K\) generates \(\mathbb H\), then \(U\) is invariant under right multiplication by every quaternion.
For every \(y\in\mathbb H\), we have \(y=u(u^{-1}y)\in U\).
Hence \(U=\mathbb H\).

By Proposition~\ref{prop:local-morse-general}, every \(G\)-invariant smooth function with a nondegenerate critical point at the origin has Morse index \(0\) or \(4\).
\end{proof}

The group \(Q_8\) generates \(\mathbb H\) as a real algebra.
For \(q\ge2\), the binary dihedral group \(\BDih{4q}=\langle e^{\pi \mathbf{i}/q},\mathbf{j}\rangle\) also generates \(\mathbb H\) as a real algebra.
Indeed, \(\mathbf{i}=(e^{\pi \mathbf{i}/q}-\cos(\pi/q))/\sin(\pi/q)\) and \(\mathbf{k}=\mathbf{i}\mathbf{j}\).
The groups \(T^*\), \(O^*\), and \(I^*\) generate \(\mathbb H\) as a real algebra because each contains a subgroup isomorphic to \(Q_8\).
Appendix~\ref{app:morse-index-case-analysis} determines the invariant subspaces for the remaining families in the completed Du Val list by computing the actions of explicit generators.

The following lemmas provide criteria used to determine the \(G\)-invariant real subspaces of \(\mathbb R^4\cong\mathbb H\).
For \(\theta\in\mathbb R\), define
\[
\operatorname{Rot}(\theta)
=
\begin{bmatrix}
\cos\theta&-\sin\theta\\
\sin\theta&\cos\theta
\end{bmatrix}.
\]

\begin{lemma}\label{lem:rot-real-eigen}
For \(\theta\in\mathbb R\), the rotation \(\operatorname{Rot}(\theta)\) has a real eigenvector if and only if \(\theta\equiv0\) or \(\pi\pmod{2\pi}\).
\end{lemma}

\begin{proof}
The discriminant of the characteristic polynomial of \(\operatorname{Rot}(\theta)\) is \(-4\sin^2\theta\).
This discriminant vanishes if and only if \(\theta\equiv0\) or \(\pi\pmod{2\pi}\).
For precisely these values, \(\operatorname{Rot}(\theta)\) has a real eigenvector.
\end{proof}

\begin{lemma}\label{lem:matrix-commute}
Let \(\alpha,\beta\in\mathbb R\).
Let \(g\in\SO(4)\) act on \(\mathbb H=\mathbb C\oplus\mathbb C\mathbf{j}\) by
\[
g(z_1+z_2\mathbf{j})=e^{\alpha \mathbf{i}}z_1+e^{\beta \mathbf{i}}z_2\mathbf{j}.
\]
Assume that neither \(\alpha\) nor \(\beta\) is congruent to \(0\) or \(\pi\) modulo \(2\pi\).
With respect to the basis \((1,\mathbf{i},\mathbf{j},\mathbf{k})\), the matrix of \(g\) is \(\operatorname{diag}\bigl(\operatorname{Rot}(\alpha),\operatorname{Rot}(\beta)\bigr)\).

The transformation \(g\) has no \(1\)- or \(3\)-dimensional real invariant subspace.
Every \(2\)-dimensional \(g\)-invariant subspace \(W\) has exactly one of the following forms.
\begin{enumerate}[label=\textup{(\roman*)}]
\item \(W=\{0\}\oplus\mathbb R^2=\operatorname{span}\{\mathbf{j},\mathbf{k}\}\).
\item \(W\) is the graph of a linear map \(T\colon\mathbb R^2\to\mathbb R^2\) satisfying \(\operatorname{Rot}(\beta)T=T\operatorname{Rot}(\alpha)\).
\end{enumerate}
\end{lemma}

\begin{proof}
By Lemma~\ref{lem:rot-real-eigen}, neither \(\operatorname{Rot}(\alpha)\) nor \(\operatorname{Rot}(\beta)\) has a real eigenvector.
Hence \(g\) has no invariant line.
Since the orthogonal complement of a \(g\)-invariant subspace is invariant, \(g\) has no invariant \(3\)-dimensional subspace.

Let \(W\subset\mathbb R^2\oplus\mathbb R^2\) be a \(2\)-dimensional \(g\)-invariant subspace.
Let \(p_1\) and \(p_2\) be the coordinate projections.
Since \(p_1(W)\) is \(\operatorname{Rot}(\alpha)\)-invariant, it is either \(0\) or \(\mathbb R^2\).
If \(p_1(W)=0\), then \(W=\{0\}\oplus\mathbb R^2\).
If \(p_1(W)=\mathbb R^2\), then \(p_1|_W\) is an isomorphism.
In this case, \(W\) is the graph of \(T=p_2\circ(p_1|_W)^{-1}\).
The subspace \(W\) is \(g\)-invariant if and only if \(\operatorname{Rot}(\beta)T=T\operatorname{Rot}(\alpha)\).
\end{proof}

\begin{lemma}\label{lem:rot-invar-basis}
Let \(\alpha,\beta\in\mathbb R\).
Let \(T\colon\mathbb R^2\to\mathbb R^2\) be a real-linear map satisfying \(\operatorname{Rot}(\beta)T=T\operatorname{Rot}(\alpha)\).
Identify \(\mathbb R^2\) with \(\mathbb C\).
Write \(T(z)=c_1z+c_2\overline z\) for unique \(c_1,c_2\in\mathbb C\).
Exactly one of the following cases occurs.

\begin{enumerate}[label=\textup{(\roman*)}]
\item
If \(\alpha\equiv\beta\pmod{2\pi}\) and \(\alpha\not\equiv0,\pi\pmod{2\pi}\), then \(c_2=0\) and
\[
T=
\begin{bmatrix}
a&-b\\
b&a
\end{bmatrix}
\]
for some \(a,b\in\mathbb R\).

\item
If \(\alpha\equiv-\beta\pmod{2\pi}\) and \(\alpha\not\equiv0,\pi\pmod{2\pi}\), then \(c_1=0\) and
\[
T=
\begin{bmatrix}
a&b\\
b&-a
\end{bmatrix}
\]
for some \(a,b\in\mathbb R\).

\item
If \((\alpha,\beta)\equiv(0,0)\) or \((\pi,\pi)\pmod{2\pi}\), then \(\operatorname{Rot}(\alpha)=\operatorname{Rot}(\beta)=\pm I_2\).
Every real-linear map \(T\colon\mathbb R^2\to\mathbb R^2\) satisfies \(\operatorname{Rot}(\beta)T=T\operatorname{Rot}(\alpha)\).

\item
In every other case, \(T=0\).
\end{enumerate}
\end{lemma}

\begin{proof}
The intertwining relation is equivalent to \(c_1(e^{\mathbf{i}\beta}-e^{\mathbf{i}\alpha})=0\) and \(c_2(e^{\mathbf{i}\beta}-e^{-\mathbf{i}\alpha})=0\).
The four cases follow from these equations.
The matrix in case~\textup{(i)} represents multiplication by \(a+b\mathbf{i}\) over the reals.
The matrix in case~\textup{(ii)} represents \(z\mapsto(a+b\mathbf{i})\overline z\) over the reals.
\end{proof}

\begin{lemma}\label{lem:unique-one-three-splitting}
Let \(G<\SO(4)\).
Suppose that \(\mathbb R^4=L\oplus W\) is an orthogonal \(G\)-invariant splitting with \(\dim L=1\) and \(\dim W=3\).
Assume that the representation of \(G\) on \(W\) is irreducible.
Then \(L\) is the unique \(G\)-invariant line in \(\mathbb R^4\).
The subspace \(W\) is the unique \(G\)-invariant \(3\)-dimensional subspace.
The splitting \(\mathbb R^4=L\oplus W\) is the unique \(G\)-invariant orthogonal splitting with one summand of dimension \(1\) and the other of dimension \(3\).
\end{lemma}

\begin{proof}
Let \(L'\subset\mathbb R^4\) be a \(G\)-invariant line.
Let \(\pi_W\colon\mathbb R^4\longrightarrow W\) be the orthogonal projection.
Since \(L\) and \(W\) are \(G\)-invariant and the action is orthogonal, \(\pi_W\) is \(G\)-equivariant.
The subspace \(\pi_W(L')\) is \(G\)-invariant and has dimension at most one.
Since \(W\) is irreducible and \(3\)-dimensional, \(\pi_W(L')=0\).
Thus \(L'\subset\ker\pi_W=L\), so \(L'=L\).

Let \(P\subset\mathbb R^4\) be a \(G\)-invariant \(3\)-dimensional subspace.
Since the action is orthogonal, \(P^\perp\) is a \(G\)-invariant line.
By the uniqueness of \(L\), we have \(P^\perp=L\).
Hence \(P=L^\perp=W\).
Thus \(\mathbb R^4=L\oplus W\) is the unique \(G\)-invariant orthogonal splitting with summands of dimensions \(1\) and \(3\).
\end{proof}

\subsection{Parameter conditions for \textup{Families} 1 and 1'}\label{subsec:parameter-cases}

The family parameters are those of Definition~\ref{def:representative-family-parameters}.
All matrices below are written in the basis \((1,\mathbf{i},\mathbf{j},\mathbf{k})\).

By Proposition~\ref{prop:generator}, the representative \(G_1\) is generated by
\begin{align*}
\gamma_1&=\operatorname{diag}\bigl(\operatorname{Rot}(\pi/m),\operatorname{Rot}(\pi/m)\bigr),\\
\gamma_2&=\operatorname{diag}\bigl(\operatorname{Rot}(-\pi/n),\operatorname{Rot}(\pi/n)\bigr),\\
\gamma_3&=\operatorname{diag}\bigl(\operatorname{Rot}(\pi/(mr)-s\pi/(nr)),\operatorname{Rot}(\pi/(mr)+s\pi/(nr))\bigr).
\end{align*}
For \textup{Family} \(1'\), write \(r=2r'\).
We use the same symbols for the generators of \(G_{1'}\).
By Proposition~\ref{prop:generator}, these generators are
\begin{align*}
\gamma_1&=\operatorname{diag}\bigl(\operatorname{Rot}(2\pi/m),\operatorname{Rot}(2\pi/m)\bigr),\\
\gamma_2&=\operatorname{diag}\bigl(\operatorname{Rot}(-2\pi/n),\operatorname{Rot}(2\pi/n)\bigr),\\
\gamma_3&=\operatorname{diag}\bigl(\operatorname{Rot}(\pi/(mr')-s\pi/(nr')),\operatorname{Rot}(\pi/(mr')+s\pi/(nr'))\bigr).
\end{align*}
Set \(\gamma_4=\operatorname{diag}(1,-1,1,-1)\).
By Proposition~\ref{prop:generator},
\[
G_{11}=\langle G_1,\gamma_4\rangle,
\qquad
G_{11'}=\langle G_{1'},\gamma_4\rangle.
\]
In each equality, the two groups have the same family parameters.

Definition~\ref{def:parameter-cases-1-1p} specifies the parameter conditions under which a representative in \textup{Family} \(1\) or \(1'\) has an invariant line or an invariant \(2\)-dimensional subspace distinct from \(\operatorname{span}\{1,\mathbf{i}\}\) and \(\operatorname{span}\{\mathbf{j},\mathbf{k}\}\).

\begin{definition}\label{def:parameter-cases-1-1p}
For \textup{Family} \(1\), let \(m\), \(n\), and \(r\) be positive integers, and let \(s\in\mathbb Z\) satisfy \(\gcd(s,r)=1\).

For \textup{Family} \(1'\), let \(m\), \(n\), and \(r'\) be positive integers, with \(m\) and \(n\) odd.
Let \(s\in\mathbb Z\) satisfy \(\gcd(s,2r')=1\).
Write \(r=2r'\).

The subfamilies are defined by the following conditions.

\begin{enumerate}[
label=\textup{Family}~\protect\ensuremath{1_{\mathrm{(\roman*)}}}.,
ref=\protect\ensuremath{1_{\mathrm{(\roman*)}}},
align=left,
labelindent=\parindent,
leftmargin=*,
labelsep=0.5em
]
\item\label{1-1}
\(m=1\) and \(r=1\).

\item\label{1-3}
\(n=1\), \(r=1\), and \((m,s)\ne(1,1)\).

\item\label{1-6}
\(m=n=1\) and \(\frac{1-s}{r}\in\mathbb Z\).
In this subfamily, set \(l_-=\frac{1-s}{r}\).

\item\label{1-7}
\(m=n=1\) and \(\frac{1+s}{r}\in\mathbb Z\).
\end{enumerate}

\begin{enumerate}[
label=\textup{Family}~\protect\ensuremath{1'_{\mathrm{(\roman*)}}}.,
ref=\protect\ensuremath{1'_{\mathrm{(\roman*)}}},
align=left,
labelindent=\parindent,
leftmargin=*,
labelsep=0.5em
]
\item\label{1'-1}
\(m=1\) and \(r'=1\).

\item\label{1'-3}
\(n=1\), \(r'=1\), and \((m,s)\ne(1,1)\).

\item\label{1'-6}
\(m=n=1\) and \(\frac{1-s}{r'}\in\mathbb Z\).
In this subfamily, set \(l'_-=\frac{1-s}{r'}\).

\item\label{1'-7}
\(m=n=1\) and \(\frac{1+s}{r'}\in\mathbb Z\).
In this subfamily, set \(l'_+=\frac{1+s}{r'}\).
\end{enumerate}
\end{definition}

\begin{remark}\label{rem:cyclic-subfamily-intersections}
The subspaces \(\operatorname{span}\{1,\mathbf{i}\}\) and \(\operatorname{span}\{\mathbf{j},\mathbf{k}\}\) are invariant for every representative in \textup{Families} \(1\) and \(1'\).
By Propositions~\ref{prop:cyclic-index-two-handle-orbits}, \ref{prop:case-family-1}, and~\ref{prop:case-family-1p}, each representative in \familyref{1-1}, \familyref{1-3}, \familyref{1'-1}, or \familyref{1'-3} has an invariant \(2\)-dimensional subspace distinct from \(\operatorname{span}\{1,\mathbf{i}\}\) and \(\operatorname{span}\{\mathbf{j},\mathbf{k}\}\).
In \familiesref{1-6}{1'-6}, the subspace \(\operatorname{span}\{1,\mathbf{i}\}\) splits into invariant lines.
In \familiesref{1-7}{1'-7}, the subspace \(\operatorname{span}\{\mathbf{j},\mathbf{k}\}\) splits into invariant lines.
A representative in \textup{Family} \(1\) admits an invariant splitting of type \(1\oplus1\oplus2\) if it belongs to \familyref{1-6} or \familyref{1-7}.
A representative in \textup{Family} \(1'\) admits an invariant splitting of type \(1\oplus1\oplus2\) if it belongs to \familyref{1'-6} or \familyref{1'-7}.

The four subfamilies of \textup{Family} \(1\) need not be disjoint, and the same holds for \textup{Family} \(1'\).
The parameter conditions for \familyref{1-6} and \familyref{1-7} hold simultaneously precisely when \(m=n=1\) and \(r\in\{1,2\}\).
The parameter conditions for \familyref{1'-6} and \familyref{1'-7} hold simultaneously precisely when \(m=n=1\) and \(r'\in\{1,2\}\).
In each family, none of the four subfamilies contains another.
\end{remark}

\begin{remark}\label{rem:four-line-splitting-criteria}
In \textup{Family} \(1\), the representation admits an invariant splitting of type \(1\oplus1\oplus1\oplus1\) precisely when it belongs to both \familyref{1-6} and \familyref{1-7}.
By Remark~\ref{rem:cyclic-subfamily-intersections}, this is equivalent to \(m=n=1\) and \(r\in\{1,2\}\).
In \textup{Family} \(1'\), the representation admits an invariant splitting of type \(1\oplus1\oplus1\oplus1\) precisely when it belongs to both \familyref{1'-6} and \familyref{1'-7}.
By Remark~\ref{rem:cyclic-subfamily-intersections}, this is equivalent to \(m=n=1\) and \(r'\in\{1,2\}\).

The chosen representatives satisfy \(G_{11}=\langle G_1,\operatorname{diag}(1,-1,1,-1)\rangle\) and \(G_{11'}=\langle G_{1'},\operatorname{diag}(1,-1,1,-1)\rangle\).
The matrix \(\operatorname{diag}(1,-1,1,-1)\) preserves \(\operatorname{span}\{1,\mathbf{i}\}\) and \(\operatorname{span}\{\mathbf{j},\mathbf{k}\}\).
The same matrix preserves each of the lines \(\operatorname{span}\{1\}\), \(\operatorname{span}\{\mathbf{i}\}\), \(\operatorname{span}\{\mathbf{j}\}\), and \(\operatorname{span}\{\mathbf{k}\}\).

In \textup{Family} \(11\), the representation admits an invariant splitting of type \(1\oplus1\oplus1\oplus1\) precisely when \(G_1\subset G_{11}\) belongs to both \familyref{1-6} and \familyref{1-7}.
In \textup{Family} \(11'\), the representation admits an invariant splitting of type \(1\oplus1\oplus1\oplus1\) precisely when \(G_{1'}\subset G_{11'}\) belongs to both \familyref{1'-6} and \familyref{1'-7}.

The restrictions \(mr,nr\ge2\) exclude \familiesref{1-1}{1-3} for \(G_1\subset G_{11}\).
The restrictions \(mr',nr'\ge2\) exclude \familiesref{1'-1}{1'-3} for \(G_{1'}\subset G_{11'}\).
Propositions~\ref{prop:case-family-1}, \ref{prop:case-family-1p}, \ref{prop:case-family-11}, and~\ref{prop:case-family-11p} establish the splitting criteria for \textup{Families} \(1\), \(1'\), \(11\), and \(11'\).
\end{remark}

\begin{remark}\label{rem:families-3-34-relation}
In \textup{Family} \(3\), \(m\ge1\) and \(n\ge2\).
In \textup{Family} \(34\), \(m\) and \(n\) are odd, with \(m\ge1\) and \(n\ge3\).
For fixed odd \(m,n\), the chosen representative \(G_{34}\) has index two in \(G_3\).
In the proofs of Propositions~\ref{prop:case-family-3} and~\ref{prop:case-family-34}, each of the first two generators of \(\widetilde G_{34}\) is the square of the corresponding generator of \(\widetilde G_3\).
The third generators of \(\widetilde G_3\) and \(\widetilde G_{34}\) agree.
Table~\ref{tab:mecchia-seppi-families} records the orders \(|G_3|=4mn\) and \(|G_{34}|=2mn\).
The representations in both families have no invariant lines or invariant \(3\)-dimensional subspaces.
By Propositions~\ref{prop:case-family-3} and~\ref{prop:case-family-34}, the representations admit an invariant splitting of type \(2\oplus2\) precisely when \(m=1\).
For all other admissible choices of parameters, the representations are irreducible.
\end{remark}

For \textup{Family} \(1'\) with \(r'>2\), \familyref{1'-6} and \familyref{1'-7} are disjoint.
A representative in both subfamilies would satisfy \(r'\mid2\).
If \(r'>2\) and the representative belongs to at least one of \familiesref{1'-6}{1'-7}, define
\[
\ell'=
\begin{cases}
l'_- & \text{in \familyref{1'-6}},\\
l'_+ & \text{in \familyref{1'-7}}.
\end{cases}
\]

\begin{lemma}\label{lem:line-summand-inheritance-11-11p}
The subfamilies are those of Definition~\ref{def:parameter-cases-1-1p}.
A representative in \textup{Family} \(11\) has an invariant line if and only if the subgroup \(G_1\subset G_{11}\) belongs to \familyref{1-6} or \familyref{1-7}.
A representative in \textup{Family} \(11'\) has an invariant line if and only if the subgroup \(G_{1'}\subset G_{11'}\) belongs to \familyref{1'-6} or \familyref{1'-7}.
\end{lemma}

\begin{proof}
The groups satisfy \(G_{11}=\langle G_1,\gamma_4\rangle\) and \(G_{11'}=\langle G_{1'},\gamma_4\rangle\).
If \(G_{11}\) has an invariant line, then that line is \(G_1\)-invariant.
By Proposition~\ref{prop:case-family-1}, \(G_1\) then belongs to \familyref{1-6} or \familyref{1-7}.
If \(G_{11'}\) has an invariant line, then that line is \(G_{1'}\)-invariant.
By Proposition~\ref{prop:case-family-1p}, \(G_{1'}\) then belongs to \familyref{1'-6} or \familyref{1'-7}.

Conversely, suppose that \(G_1\) belongs to \familyref{1-6} or \familyref{1-7}.
If \(G_1\) belongs to \familyref{1-6}, then \(\operatorname{span}\{1\}\) is \(G_1\)-invariant.
If \(G_1\) belongs to \familyref{1-7}, then \(\operatorname{span}\{\mathbf{j}\}\) is \(G_1\)-invariant.
Since \(\gamma_4\) preserves both lines, \(G_{11}\) has an invariant line in either case.

Suppose that \(G_{1'}\) belongs to \familyref{1'-6} or \familyref{1'-7}.
If \(G_{1'}\) belongs to \familyref{1'-6}, then \(\operatorname{span}\{1\}\) is \(G_{1'}\)-invariant.
If \(G_{1'}\) belongs to \familyref{1'-7}, then \(\operatorname{span}\{\mathbf{j}\}\) is \(G_{1'}\)-invariant.
Since \(\gamma_4\) preserves both lines, \(G_{11'}\) has an invariant line in either case.
\end{proof}

\subsection{Classification of Morse indices}\label{subsec:morse-index-theorem}

\begin{theorem}\label{thm:4d-morse-index-classification}
Let \(G<\SO(4)\) be one of the representatives fixed in Section~\ref{subsec:so4-finite-subgroups}.
Assume that the parameters satisfy the conditions in Definition~\ref{def:representative-family-parameters}.
The possible Morse indices of \(G\)-invariant smooth functions \(F\colon\mathbb R^4\to\mathbb R\) with a nondegenerate critical point at the origin are as follows.

\begin{enumerate}[label=\textup{(\roman*)}]
\item
The following representations are irreducible.
The possible Morse indices are \(0\) and \(4\).
\begin{flushleft}
\upshape
\textup{Families} \(5\)--\(9\); \textup{Families} \(14\)--\(20\); \textup{Families} \(22\)--\(25\);\\
\textup{Families} \(27\)--\(30\);\\
\textup{Families} \(2\), \(4\), \(10\), and \(12\); \textup{Family} \(13\); \textup{Families} \(32\) and \(32'\);\\
\textup{Families} \(33\) and \(33'\);\\
\textup{Family} \(3\) with \(m,n\ge2\);\\
\textup{Family} \(34\) with \(m,n\ge3\).
\end{flushleft}

\item
The following representations are direct sums of two irreducible real \(2\)-dimensional modules.
The possible Morse indices are \(0\), \(2\), and \(4\).
\begin{flushleft}
\upshape
\textup{Family} \(1\) when \(G\) belongs to neither \familyref{1-6} nor \familyref{1-7};\\
\textup{Family} \(1'\) when \(G\) belongs to neither \familyref{1'-6} nor \familyref{1'-7};\\
\textup{Family} \(11\) when \(G_1\subset G_{11}\) belongs to neither \familyref{1-6} nor \familyref{1-7};\\
\textup{Family} \(11'\) when \(G_{1'}\subset G_{11'}\) belongs to neither \familyref{1'-6} nor \familyref{1'-7};\\
\textup{Family} \(3\) with \(m=1\) and \(n\ge2\);\\
\textup{Family} \(34\) with \(m=1\) and odd \(n\ge3\).
\end{flushleft}

\item
The representations in \textup{Families} \(21\), \(21'\), \(26\), \(26'\), \(26''\), \(31\), and \(31'\) admit an invariant splitting of type \(1\oplus3\).
The \(3\)-dimensional summand is irreducible.
The possible Morse indices are \(0\), \(1\), \(3\), and \(4\).

\item
The following representations admit an invariant splitting of type \(1\oplus1\oplus2\).
The possible Morse indices are \(0\), \(1\), \(2\), \(3\), and \(4\).
\begin{flushleft}
\upshape
\textup{Family} \(1\) when \(G\) belongs to \familyref{1-6} or \familyref{1-7};\\
\textup{Family} \(1'\) when \(G\) belongs to \familyref{1'-6} or \familyref{1'-7};\\
\textup{Family} \(11\) when \(G_1\subset G_{11}\) belongs to \familyref{1-6} or \familyref{1-7};\\
\textup{Family} \(11'\) when \(G_{1'}\subset G_{11'}\) belongs to \familyref{1'-6} or \familyref{1'-7}.
\end{flushleft}
\end{enumerate}
\end{theorem}

\begin{proposition}\label{prop:criterion-q8}
Let \(G<\SO(4)\) be a chosen representative in one of the following families.
\begin{itemize}[label={}]
\item \textup{Families \(2\) and \(4\).}
\item \textup{Families \(5\)--\(10\) and \(12\)--\(20\).}
\item \textup{Families \(22\)--\(25\) and \(27\)--\(30\).}
\end{itemize}
Then the real \(G\)-module \(\mathbb R^4\cong\mathbb H\) is irreducible.
Every \(G\)-invariant smooth function \(F\colon\mathbb R^4\to\mathbb R\) with a nondegenerate critical point at the origin has Morse index \(0\) or \(4\).
\end{proposition}

\begin{proof}
For each family in the proposition, Table~\ref{tab:mecchia-seppi-families} identifies either \(L_K\) or \(R_K\) with one of the groups \(\BDih{4q}\) for \(q\ge2\), \(T^*\), \(O^*\), or \(I^*\).
Each group generates \(\mathbb H\) as a real algebra.
By Proposition~\ref{prop:irred-local-gps}, the real \(G\)-module \(\mathbb H\) is irreducible.
By Proposition~\ref{prop:local-morse-general}, every \(G\)-invariant smooth function with a nondegenerate critical point at the origin has Morse index \(0\) or \(4\).
The quadratic forms \(Q_{0,\mathbb R^4}\) and \(Q_{\mathbb R^4,0}\) realize these two indices.
\end{proof}

\begin{definition}\label{def:block-graph-plane}
Using the basis \((1,\mathbf{i},\mathbf{j},\mathbf{k})\), identify \(\mathbb H\) with \(\mathbb R^2\oplus\mathbb R^2\).
For \(\theta\in\mathbb R\), set
\[
M_\theta
=
\begin{bmatrix}
\cos\theta&\sin\theta\\
\sin\theta&-\cos\theta
\end{bmatrix},
\qquad
W_\theta
=
\{(v,M_\theta v)\mid v\in\mathbb R^2\}.
\]
Write \(W_\theta^\perp\) for the orthogonal complement of \(W_\theta\).
\end{definition}

\begin{proposition}\label{prop:criterion-line-summand}
Let \(G<\SO(4)\) be a chosen representative in one of the following families.
\begin{itemize}[label={}]
\item \textup{Family \(1\) when \(G\) belongs to \familyref{1-6} or \familyref{1-7}.}
\item \textup{Family \(1'\) when \(G\) belongs to \familyref{1'-6} or \familyref{1'-7}.}
\item \textup{Family \(11\) when the subgroup \(G_1\subset G_{11}\) belongs to \familyref{1-6} or \familyref{1-7}.}
\item \textup{Family \(11'\) when the subgroup \(G_{1'}\subset G_{11'}\) belongs to \familyref{1'-6} or \familyref{1'-7}.}
\end{itemize}
Assume that the parameters of \(G\) satisfy the conditions in Definition~\ref{def:representative-family-parameters}.
Then the real \(G\)-module \(\mathbb R^4\) admits an invariant splitting
\[
\mathbb R^4=L_1\oplus L_2\oplus W,
\qquad
\dim L_1=\dim L_2=1,
\qquad
\dim W=2.
\]
Each \(k\in\{0,1,2,3,4\}\) occurs as the Morse index at the origin of a \(G\)-invariant nondegenerate quadratic form on \(\mathbb R^4\).
\end{proposition}

\begin{proof}
Suppose first that \(G\) is in \textup{Family} \(1\).
If \(G\) belongs to \familyref{1-6}, then \(\operatorname{span}\{1,\mathbf{i}\}\) splits into two invariant lines and \(\operatorname{span}\{\mathbf{j},\mathbf{k}\}\) is invariant.
If \(G\) belongs to \familyref{1-7}, then \(\operatorname{span}\{\mathbf{j},\mathbf{k}\}\) splits into two invariant lines and \(\operatorname{span}\{1,\mathbf{i}\}\) is invariant.
By Proposition~\ref{prop:case-family-1}, the representation admits an invariant splitting of type \(1\oplus1\oplus2\).

By Proposition~\ref{prop:case-family-1p}, the representation in \textup{Family} \(1'\) admits an invariant splitting of type \(1\oplus1\oplus2\) when \(G\) belongs to \familyref{1'-6} or \familyref{1'-7}.
Propositions~\ref{prop:case-family-11} and~\ref{prop:case-family-11p} establish the same splitting type for \textup{Families} \(11\) and \(11'\), respectively.

The sums of summands in \(L_1\oplus L_2\oplus W\) include invariant subspaces of every dimension from \(0\) to \(4\).
By Lemma~\ref{lem:main-hessian-splitting-principle}, each dimension occurs as a Morse index.
\end{proof}

\begin{proposition}\label{prop:criterion-polyhedral-one-three}
Let \(G<\SO(4)\) be a chosen representative in one of \textup{Families} \(21\), \(21'\), \(26\), \(26'\), \(26''\), \(31\), and \(31'\).
The real \(G\)-module \(\mathbb H\) admits the following invariant splittings.

\begin{enumerate}[label=\textup{(\roman*)}]
\item
For \textup{Families} \(21\), \(21'\), \(26\), \(26'\), \(31\), and \(31'\),
\[
\mathbb H
=
\operatorname{span}\{1\}
\oplus
\operatorname{span}\{\mathbf{i},\mathbf{j},\mathbf{k}\}.
\]

\item
For \textup{Family} \(26''\),
\[
\mathbb H
=
\operatorname{span}\{\mathbf{i}\}
\oplus
\operatorname{span}\{1,\mathbf{j},\mathbf{k}\}.
\]
\end{enumerate}

In each case, the \(3\)-dimensional summand is irreducible.
Every \(G\)-invariant smooth function \(F\colon\mathbb R^4\to\mathbb R\) with a nondegenerate critical point at the origin has Morse index \(0\), \(1\), \(3\), or \(4\).
\end{proposition}

\begin{proof}
By Propositions~\ref{prop:case-family-21-21p}, \ref{prop:case-family-26}, and~\ref{prop:case-family-31-31p}, each splitting in the proposition is \(G\)-invariant, with an irreducible \(3\)-dimensional summand.
By Proposition~\ref{prop:local-morse-general}, every \(G\)-invariant smooth function with a nondegenerate critical point at the origin has Morse index \(0\), \(1\), \(3\), or \(4\).

Write the invariant splitting as \(\mathbb H=L\oplus W\), with \(\dim L=1\) and \(\dim W=3\).
By Lemma~\ref{lem:main-hessian-splitting-principle}, \(Q_{V,V^\perp}\) realizes Morse index \(\dim V\) for each \(V\in\{0,L,W,\mathbb H\}\).
\end{proof}

\begin{corollary}\label{cor:unique-index-one-three-TOI}
Let \(G<\SO(4)\) be a chosen representative in \textup{Family} \(21\), \(21'\), \(26\), \(26'\), \(26''\), \(31\), or \(31'\).
Let \(\rho\colon G\hookrightarrow\SO(4)\) be the inclusion.
There is a unique \(G\)-invariant line \(L\) and a unique \(G\)-invariant \(3\)-dimensional subspace \(W\) in \(\mathbb H\).
These subspaces are as follows.

\begin{enumerate}[label=\textup{(\roman*)}]
\item
For \textup{Families} \(21\), \(21'\), \(26\), \(26'\), \(31\), and \(31'\),
\[
L=\operatorname{span}\{1\},
\qquad
W=\operatorname{span}\{\mathbf{i},\mathbf{j},\mathbf{k}\}.
\]

\item
For \textup{Family} \(26''\),
\[
L=\operatorname{span}\{\mathbf{i}\},
\qquad
W=\operatorname{span}\{1,\mathbf{j},\mathbf{k}\}.
\]
\end{enumerate}

The representation \(\rho\) determines a unique isomorphism class of ordered handles for each of the indices \(1\) and \(3\).
\end{corollary}

\begin{proof}
By Proposition~\ref{prop:criterion-polyhedral-one-three} and Lemma~\ref{lem:unique-one-three-splitting}, \(L\) is the unique \(G\)-invariant line and \(W\) is the unique \(G\)-invariant \(3\)-dimensional subspace.
Hence the orbits of \(L\) and \(W\) under \(N^+_\rho\) are singletons.
By Proposition~\ref{prop:ordered-handle-orbits}, \(\rho\) determines a unique isomorphism class of ordered handles for each of the indices \(1\) and \(3\).
\end{proof}

\begin{proof}[\textup{Proof of Theorem~\ref{thm:4d-morse-index-classification}}]
For the parameter ranges in part~\textup{(i)}, the representations are irreducible by Proposition~\ref{prop:criterion-q8} and Propositions~\ref{prop:case-family-3}, \ref{prop:case-family-33-33p}, \ref{prop:case-family-34}, and~\ref{prop:case-family-32-32p}.
By Proposition~\ref{prop:local-morse-general}, the possible Morse indices are \(0\) and \(4\).
The quadratic forms \(Q_{0,\mathbb R^4}\) and \(Q_{\mathbb R^4,0}\) realize these indices.

By Propositions~\ref{prop:case-family-1}, \ref{prop:case-family-1p}, \ref{prop:case-family-3}, \ref{prop:case-family-34}, \ref{prop:case-family-11}, and~\ref{prop:case-family-11p}, each representation in part~\textup{(ii)} has an invariant \(2\)-dimensional subspace \(V\).
These propositions, together with Lemma~\ref{lem:line-summand-inheritance-11-11p}, show that \(G\) has no invariant line.
The orthogonal complement \(V^\perp\) is also \(G\)-invariant.
Neither \(V\) nor \(V^\perp\) contains a \(G\)-invariant line.
Hence both summands are irreducible.
By Proposition~\ref{prop:local-morse-general}, the possible Morse indices are \(0\), \(2\), and \(4\).
The quadratic forms \(Q_{0,\mathbb R^4}\), \(Q_{V,V^\perp}\), and \(Q_{\mathbb R^4,0}\) realize these indices.

Parts~\textup{(iii)} and~\textup{(iv)} follow from Propositions~\ref{prop:criterion-polyhedral-one-three} and~\ref{prop:criterion-line-summand}, respectively.
\end{proof}

\section{Classification of four-dimensional ordered handles}\label{sec:four-dimensional-handle-classification}

In this section, we classify the ordered handles determined by the invariant splittings in Section~\ref{sec:morse-index-classification}.
By Proposition~\ref{prop:ordered-handle-orbits}, it suffices to determine the normalizer orbits of invariant subspaces for each chosen representative.
Lemma~\ref{lem:coordinate-factor-swap} determines when an element of the normalizer interchanges \(\operatorname{span}\{1,\mathbf{i}\}\) and \(\operatorname{span}\{\mathbf{j},\mathbf{k}\}\).
Section~\ref{subsec:cyclic-plane-orbits} classifies the orbits of invariant \(2\)-dimensional subspaces in \textup{Families} \(1\) and \(1'\).
Section~\ref{subsec:invariant-line-orbits} classifies the orbits of invariant lines in \textup{Families} \(1\), \(1'\), \(11\), and \(11'\).
Section~\ref{subsec:remaining-plane-orbits} determines the orbits of invariant \(2\)-dimensional subspaces in \textup{Families} \(11\), \(11'\), \(3\), and \(34\).
The proof of Theorem B concludes the section.

We use the generators \(\gamma_1,\gamma_2,\gamma_3\) and \(\gamma_4\) from Section~\ref{subsec:parameter-cases}.

\begin{lemma}\label{lem:coordinate-factor-swap}
Let \(G<\SO(4)\) be a chosen representative in one of \textup{Families} \(1\), \(11\), \(1'\), or \(11'\).
Whenever they occur, \(m,n,r,s,r'\) are the parameters of \(G\) in Definition~\ref{def:representative-family-parameters}.
Let \(\rho\colon G\hookrightarrow\SO(4)\) be the inclusion.
For \textup{Families} \(1\) and \(11\), the ordered handles
\[
\mathbf H_\rho(\operatorname{span}\{1,\mathbf{i}\})
\qquad\text{and}\qquad
\mathbf H_\rho(\operatorname{span}\{\mathbf{j},\mathbf{k}\})
\]
are isomorphic precisely when \(r=1\) or \(r=2\).
For \textup{Families} \(1'\) and \(11'\), these ordered handles are isomorphic precisely when \(r'=1\).
\end{lemma}

\begin{proof}
By Proposition~\ref{prop:ordered-handle-orbits}, it suffices to determine when an element of \(N^+_\rho\) interchanges \(\operatorname{span}\{1,\mathbf{i}\}\) and \(\operatorname{span}\{\mathbf{j},\mathbf{k}\}\).

Set \(S^1=S^3\cap\mathbb C\) and
\[
K=\SO(\operatorname{span}\{1,\mathbf{i}\})
\times
\SO(\operatorname{span}\{\mathbf{j},\mathbf{k}\}).
\]
For \(p,q\in S^1\),
\[
\Phi(p,q)(z_1+z_2\mathbf{j})
=
(pq^{-1})z_1+(pq)z_2\mathbf{j}.
\]
Let \(a,b\in S^1\).
Choose \(p\in S^1\) with \(p^2=ab\), and set \(q=p^{-1}b\).
Then \(pq^{-1}=a\) and \(pq=b\).
Thus \(\Phi(S^1\times S^1)=K\).

Define \(\beta\colon\mathbb H\to\mathbb H\) by
\[
\beta(z_1+z_2\mathbf{j})=z_2+z_1\mathbf{j},
\qquad
z_1,z_2\in\mathbb C.
\]
Then \(\beta\in\SO(4)\).
It satisfies \(\beta(\operatorname{span}\{1,\mathbf{i}\}) =\operatorname{span}\{\mathbf{j},\mathbf{k}\}\).
For \(p,q\in S^1\),
\[
\beta\Phi(p,q)\beta^{-1}
=
\Phi(p,q^{-1}).
\]

With respect to the splitting
\[
\mathbb H
=
\operatorname{span}\{1,\mathbf{i}\}
\oplus
\operatorname{span}\{\mathbf{j},\mathbf{k}\},
\]
every \(T\in\SO(4)\) interchanging the two summands has the form
\[
T=
\begin{pmatrix}
0&A\\
B&0
\end{pmatrix},
\qquad
A,B\in\Orth(2),
\qquad
\det A=\det B=\varepsilon\in\{-1,1\}.
\]
Conjugation by an element of \(\Orth(2)\) with determinant \(\varepsilon\) sends every rotation in \(\SO(2)\) to its \(\varepsilon\)-th power.
Hence \(T\Phi(p,q)T^{-1}=\Phi(p^\varepsilon,q^{-\varepsilon})\) for all \(p,q\in S^1\).
If \(\varepsilon=1\), conjugation by \(T\) agrees on \(K\) with conjugation by \(\beta\).
If \(\varepsilon=-1\), then
\[
T\Phi(p,q)T^{-1}
=
\bigl(\beta\Phi(p,q)\beta^{-1}\bigr)^{-1}.
\]
Since \(H^{-1}=H\) for every subgroup \(H<K\), we have \(THT^{-1}=\beta H\beta^{-1}\).
The map \(T\) normalizes \(H\) if and only if \(\beta\) normalizes \(H\).

For \textup{Family} \(1\), let \(L,L_K,R,R_K\) be the subgroups in the description of \(\widetilde G_1\) from Proposition~\ref{prop:quaternionic-cover-fiber-product}.
Choose generators \(\overline u\) of \(L/L_K\) and \(\overline v\) of \(R/R_K\) such that \(\varphi_{\widetilde G_1}(\overline u)=\overline v^{\,s}\).
For \((p,q)\in L\times R\), write \(pL_K=\overline u^{\,a}\) and \(qR_K=\overline v^{\,b}\).
The pair \((p,q)\) lies in \(\widetilde G_1\) if and only if \(b\equiv sa\pmod r\).
Conjugation by \(\beta\) corresponds to \((p,q)\mapsto(p,q^{-1})\).
This map sends \((a,b)\) to \((a,-b)\).
Since \(\ker\Phi\subset\widetilde G_1\), Proposition~\ref{prop:quaternionic-cover-fiber-product} implies that \(\beta G_1\beta^{-1}=G_1\) if and only if this map preserves \(\widetilde G_1\).
This condition is equivalent to \(2s\equiv0\pmod r\).
Since \(\gcd(r,s)=1\) by Definition~\ref{def:representative-family-parameters}, this is equivalent to \(r\in\{1,2\}\).

For \textup{Family} \(1'\), let \(L,L_K,R,R_K\) be the subgroups in the description of \(\widetilde G_{1'}\).
The quotient groups \(L/L_K\) and \(R/R_K\) have order \(2r'\).
Choose generators \(\overline u\) of \(L/L_K\) and \(\overline v\) of \(R/R_K\) such that \(\varphi_{\widetilde G_{1'}}(\overline u)=\overline v^{\,s}\).
For \((p,q)\in L\times R\), write \(pL_K=\overline u^{\,a}\) and \(qR_K=\overline v^{\,b}\), where \(a,b\in\mathbb Z/(2r')\mathbb Z\).
The pair \((p,q)\) lies in \(\widetilde G_{1'}\) if and only if \(b\equiv sa\pmod{2r'}\).
Conjugation by \(\beta\) preserves \(G_{1'}\) if and only if \(2s\equiv0\pmod{2r'}\).
Since \(\gcd(s,2r')=1\) by Definition~\ref{def:representative-family-parameters}, this is equivalent to \(r'=1\).

By the generator formulas in Section~\ref{subsec:parameter-cases}, \(G_{11}=\langle G_1,\gamma_4\rangle\) and \(G_{11'}=\langle G_{1'},\gamma_4\rangle\).
For \(g=\Phi(p,q)\in G_1\) or \(G_{1'}\), with \(p,q\in S^1\), conjugation gives \(\gamma_4g\gamma_4^{-1}=\Phi(p^{-1},q^{-1})=g^{-1}\).
Thus \(\gamma_4\) normalizes both \(G_1\) and \(G_{1'}\).
The groups \(G_1\) and \(G_{1'}\) lie in \(K\).
The restriction of \(\gamma_4\) to each of \(\operatorname{span}\{1,\mathbf{i}\}\) and \(\operatorname{span}\{\mathbf{j},\mathbf{k}\}\) has determinant \(-1\).
Therefore
\[
G_1=G_{11}\cap K,
\qquad
G_{1'}=G_{11'}\cap K.
\]

Let \(T\in\SO(4)\) interchange \(\operatorname{span}\{1,\mathbf{i}\}\) and \(\operatorname{span}\{\mathbf{j},\mathbf{k}\}\).
Then \(TKT^{-1}=K\).
If \(T\) normalizes \(G_{11}\), then
\[
TG_1T^{-1}
=
T(G_{11}\cap K)T^{-1}
=
G_{11}\cap K
=
G_1.
\]
The calculation for \textup{Family} \(1\) therefore implies \(r\in\{1,2\}\).
If \(T\) normalizes \(G_{11'}\), then \(TG_{1'}T^{-1}=G_{1'}\).
The calculation for \textup{Family} \(1'\) therefore implies \(r'=1\).

Conversely, if \(r\in\{1,2\}\), then \(\beta G_1\beta^{-1}=G_1\).
If \(r'=1\), then \(\beta G_{1'}\beta^{-1}=G_{1'}\).
Since \(\beta\gamma_4=\gamma_4\beta\), we have \(\beta G_{11}\beta^{-1}=G_{11}\) when \(r\in\{1,2\}\).
The same identity implies \(\beta G_{11'}\beta^{-1}=G_{11'}\) when \(r'=1\).
The map \(\beta\) interchanges \(\operatorname{span}\{1,\mathbf{i}\}\) and \(\operatorname{span}\{\mathbf{j},\mathbf{k}\}\).

Conjugation by \(\beta\) replaces \(s\) by \(-s\).
It exchanges the parameter ranges defining \familyref{1-6} and \familyref{1-7}.
It also exchanges the parameter ranges defining \familyref{1'-6} and \familyref{1'-7}.
\end{proof}

\subsection{\textup{Families} \texorpdfstring{\(1\) and \(1'\)}{1 and 1'}}\label{subsec:cyclic-plane-orbits}

For \(G<\SO(4)\), the centralizer of \(G\) in \(\SO(4)\) is
\[
C_{\SO(4)}(G)
=
\{S\in\SO(4)\mid Sg=gS\text{ for every }g\in G\}.
\]

\begin{lemma}\label{lem:two-block-invariant-plane-orbits}
Let \(G<\SO(4)\) be finite.
Let \(\mathbb R^4=E\oplus F\) be an orthogonal \(G\)-invariant splitting with \(\dim E=\dim F=2\).
Let \(P\subset\mathbb R^4\) be a \(G\)-invariant \(2\)-dimensional subspace.
After interchanging \(E\) and \(F\) if necessary, one of the following holds.
\begin{enumerate}[label=\textup{(\roman*)}]
\item
The subspace \(P\) is the graph of a \(G\)-equivariant linear map \(T\colon E\to F\).

\item
There are \(G\)-invariant lines \(\ell_E\subset E\) and \(\ell_F\subset F\) such that \(P=\ell_E\oplus\ell_F\).
\end{enumerate}

Suppose that \(P=\operatorname{graph}(T)\).
Then there exists \(S\in C_{\SO(4)}(G)\) such that \(S(P)=E\) and \(S(P^\perp)=F\).
In particular, \(P\) and \(E\) lie in the same \(N_{\SO(4)}(G)\)-orbit.
\end{lemma}

\begin{proof}
Let \(\pi_E\colon\mathbb R^4\to E\) and \(\pi_F\colon\mathbb R^4\to F\) be the orthogonal projections.
Their restrictions to \(P\) are \(G\)-equivariant.

Suppose that one of \(\pi_E|_P\) and \(\pi_F|_P\) has rank \(2\).
After interchanging \(E\) and \(F\) if necessary, assume that \(\pi_E|_P\colon P\to E\) is an isomorphism.
Then \(P=\operatorname{graph}(T)\), where \(T=\pi_F\circ(\pi_E|_P)^{-1}\colon E\to F\).
The map \(T\) is \(G\)-equivariant.

Suppose that neither restriction has rank \(2\).
Neither restriction has rank \(0\).
If \(\pi_E|_P=0\), then \(P=F\), and hence \(\pi_F|_P\) has rank \(2\).
The same argument applies with \(E\) and \(F\) interchanged.
Both restrictions therefore have rank \(1\).
The images \(\ell_E=\pi_E(P)\) and \(\ell_F=\pi_F(P)\) are \(G\)-invariant lines.
Every vector in \(P\) has its \(E\)-component in \(\ell_E\) and its \(F\)-component in \(\ell_F\).
Thus \(P\subset\ell_E\oplus\ell_F\).
Since both subspaces have dimension \(2\), we have \(P=\ell_E\oplus\ell_F\).

Suppose that \(P=\operatorname{graph}(T)\) for a \(G\)-equivariant linear map \(T\colon E\to F\).
With respect to the orthogonal splitting \(\mathbb R^4=E\oplus F\), define
\[
K=
\begin{pmatrix}
0&-T^*\\
T&0
\end{pmatrix},
\qquad
B=I+K.
\]
Since \(T\) is \(G\)-equivariant and the action of \(G\) is orthogonal, \(T^*\) is \(G\)-equivariant.
The maps \(K\) and \(B\) commute with \(G\).
Moreover, \(K^*=-K\).

For every \(t\in\mathbb R\), the map \(I+tK\) is invertible.
If \((I+tK)v=0\), then \(0=\langle (I+tK)v,v\rangle=\lVert v\rVert^2\), since \(\langle Kv,v\rangle=0\).
Thus \(v=0\).
It follows that \(\det(I+tK)\ne0\) for every \(t\in\mathbb R\).
Since \(\det I=1\), continuity gives \(\det B>0\).

Write the polar decomposition as \(B=UR\), where \(U\) is orthogonal and \(R=(B^*B)^{1/2}\) is positive definite.
We have
\[
B^*B
=
I-K^2
=
\begin{pmatrix}
I_E+T^*T&0\\
0&I_F+TT^*
\end{pmatrix}.
\]
The positive square root \(R\) preserves \(E\) and \(F\).
Since it commutes with \(G\), the map \(U=BR^{-1}\) also commutes with \(G\).

The map \(B\) satisfies \(B(E)=\operatorname{graph}(T)=P\) and \(B(F)=\{(-T^*y,y)\mid y\in F\}=P^\perp\).
Since \(R\) preserves \(E\) and \(F\), we have \(U(E)=P\) and \(U(F)=P^\perp\).
Since \(R\) is positive definite, \(\det R>0\).
Together with \(\det B>0\), this implies \(\det U>0\).
Thus \(U\in\SO(4)\).

Set \(S=U^{-1}\).
The map \(S\) commutes with \(G\) and satisfies \(S(P)=E\) and \(S(P^\perp)=F\).
Therefore \(S\in C_{\SO(4)}(G)\).
\end{proof}

\begin{proposition}\label{prop:cyclic-index-two-handle-orbits}
Let \(G<\SO(4)\) be a chosen representative in \textup{Family} \(1\) or \(1'\).
Whenever they occur, \(m,n,r,s,r'\) are the parameters of \(G\) in Definition~\ref{def:representative-family-parameters}.
Let \(\rho\colon G\hookrightarrow\SO(4)\) be the inclusion.

\begin{enumerate}[label=\textup{(\roman*)}]
\item
For \textup{Family} \(1\), the \(2\)-dimensional subspaces \(\operatorname{span}\{1,\mathbf{i}\}\) and \(\operatorname{span}\{\mathbf{j},\mathbf{k}\}\) lie in the same \(N^+_\rho\)-orbit precisely when \(r=1\) or \(r=2\).
For \(r>2\), they represent two distinct orbits.
Orthogonal complementation interchanges these two orbits.
When \(m=n=1\) and \(r=2\), there is one additional orbit represented by \(\ell_1\oplus\ell_2\), where \(\ell_1\subset \operatorname{span}\{1,\mathbf{i}\}\) and \(\ell_2\subset \operatorname{span}\{\mathbf{j},\mathbf{k}\}\) are invariant lines.
This orbit is distinct from the orbit of \(\operatorname{span}\{1,\mathbf{i}\}\).
Orthogonal complementation leaves this orbit invariant.

\item
For \textup{Family} \(1'\), the \(2\)-dimensional subspaces \(\operatorname{span}\{1,\mathbf{i}\}\) and \(\operatorname{span}\{\mathbf{j},\mathbf{k}\}\) lie in the same \(N^+_\rho\)-orbit precisely when \(r'=1\).
For \(r'>1\), they represent two distinct orbits.
Orthogonal complementation interchanges these two orbits.
When \(m=n=1\) and \(r'=2\), there is one additional orbit represented by \(\ell_1\oplus\ell_2\), where \(\ell_1\subset \operatorname{span}\{1,\mathbf{i}\}\) and \(\ell_2\subset \operatorname{span}\{\mathbf{j},\mathbf{k}\}\) are invariant lines.
This orbit is distinct from the orbits of \(\operatorname{span}\{1,\mathbf{i}\}\) and \(\operatorname{span}\{\mathbf{j},\mathbf{k}\}\).
Orthogonal complementation leaves this orbit invariant.
\end{enumerate}

These are all the \(N^+_\rho\)-orbits of invariant \(2\)-dimensional subspaces.
\end{proposition}

\begin{proof}
By Lemma~\ref{lem:two-block-invariant-plane-orbits}, every invariant \(2\)-dimensional subspace has one of two forms, after interchanging \(\operatorname{span}\{1,\mathbf{i}\}\) and \(\operatorname{span}\{\mathbf{j},\mathbf{k}\}\) if necessary.
It is either the graph of a \(G\)-equivariant linear map or a direct sum of invariant lines in the two subspaces.

Identify \(\operatorname{span}\{1,\mathbf{i}\}\) with \(\mathbb C\) and \(\operatorname{span}\{\mathbf{j},\mathbf{k}\}\) with \(\mathbb C\mathbf{j}\).
Every real-linear map \(T\colon \operatorname{span}\{1,\mathbf{i}\}\to \operatorname{span}\{\mathbf{j},\mathbf{k}\}\) has the form
\[
T(z)=(uz+v\overline z)\mathbf{j},
\qquad u,v\in\mathbb C.
\]
Suppose that the action is nonscalar.
By Lemma~\ref{lem:rot-invar-basis} and the generator formulas in Section~\ref{subsec:parameter-cases}, a \(G\)-equivariant map \(T\) with \(u\ne0\) exists precisely under the conditions of \familyref{1-3} for \textup{Family} \(1\) and \familyref{1'-3} for \textup{Family} \(1'\).
Similarly, a \(G\)-equivariant map \(T\) with \(v\ne0\) exists precisely under the conditions of \familyref{1-1} for \textup{Family} \(1\) and \familyref{1'-1} for \textup{Family} \(1'\).
By Lemma~\ref{lem:two-block-invariant-plane-orbits}, every invariant \(2\)-dimensional subspace that is the graph of such a map lies in the \(N^+_\rho\)-orbit of \(\operatorname{span}\{1,\mathbf{i}\}\).
For a scalar action, \(C_{\SO(4)}(G)=\SO(4)\), and \(\SO(4)\) acts transitively on \(\operatorname{Gr}_2(\mathbb R^4)\).

Suppose that the action in \textup{Family} \(1\) is nonscalar.
By Proposition~\ref{prop:case-family-1}, an invariant subspace \(P=\ell_1\oplus\ell_2\), with \(\ell_1\subset \operatorname{span}\{1,\mathbf{i}\}\) and \(\ell_2\subset \operatorname{span}\{\mathbf{j},\mathbf{k}\}\) invariant lines, exists precisely when \(m=n=1\) and \(r=2\).
For a nonscalar action in \textup{Family} \(1'\), such a subspace exists precisely when \(m=n=1\) and \(r'=2\), by Proposition~\ref{prop:case-family-1p}.
In both cases, every element of \(G\) restricts to \(\pm I\) on each of \(\operatorname{span}\{1,\mathbf{i}\}\) and \(\operatorname{span}\{\mathbf{j},\mathbf{k}\}\).
The group \(\SO(\operatorname{span}\{1,\mathbf{i}\})\times\SO(\operatorname{span}\{\mathbf{j},\mathbf{k}\})\) centralizes \(G\) and acts transitively on the subspaces \(P=\ell_1\oplus\ell_2\).

The representations of \(G\) on \(\operatorname{span}\{1,\mathbf{i}\}\) and \(\operatorname{span}\{\mathbf{j},\mathbf{k}\}\) are direct sums of two copies of real characters \(\chi_1\) and \(\chi_2\), respectively.
The representation of \(G\) on \(P=\ell_1\oplus\ell_2\) is the direct sum of \(\chi_1\) and \(\chi_2\), with \(\chi_1\ne\chi_2\).

For \(S\in N^+_\rho\), let \(\alpha_S\in\operatorname{Aut}(G)\) be the automorphism induced by conjugation.
The representation of \(G\) on \(S(P)\) is the direct sum of \(\chi_1\circ\alpha_S^{-1}\) and \(\chi_2\circ\alpha_S^{-1}\).
These characters are distinct.
On each of \(\operatorname{span}\{1,\mathbf{i}\}\) and \(\operatorname{span}\{\mathbf{j},\mathbf{k}\}\), the representation of \(G\) is the direct sum of two copies of a single character.
Therefore \(S(P)\) is neither \(\operatorname{span}\{1,\mathbf{i}\}\) nor \(\operatorname{span}\{\mathbf{j},\mathbf{k}\}\).

By Lemma~\ref{lem:coordinate-factor-swap}, \(\operatorname{span}\{1,\mathbf{i}\}\) and \(\operatorname{span}\{\mathbf{j},\mathbf{k}\}\) lie in the same \(N^+_\rho\)-orbit precisely when \(r\in\{1,2\}\) in \textup{Family} \(1\) and when \(r'=1\) in \textup{Family} \(1'\).

Let \(R_1\in\SO(\operatorname{span}\{1,\mathbf{i}\})\) and \(R_2\in\SO(\operatorname{span}\{\mathbf{j},\mathbf{k}\})\) be the rotations through \(\pi/2\).
For \(P=\ell_1\oplus\ell_2\), we have \(R_1\oplus R_2\in C_{\SO(4)}(G)\) and \((R_1\oplus R_2)(P)=P^\perp\).
Orthogonal complementation therefore leaves the \(N^+_\rho\)-orbit of \(P\) invariant.
\end{proof}

\subsection{Invariant lines}\label{subsec:invariant-line-orbits}

\begin{proposition}\label{prop:line-handle-orbits}
Let \(G<\SO(4)\) be a chosen representative satisfying the hypotheses of Proposition~\ref{prop:criterion-line-summand}.
Whenever they occur, \(m,n,r,s,r'\) are the parameters of \(G\) in Definition~\ref{def:representative-family-parameters}.
Let \(\rho\colon G\hookrightarrow\SO(4)\) be the inclusion.
The \(N^+_\rho\)-orbits of invariant lines are as follows.

\begin{enumerate}[label=\textup{(\roman*)}]
\item
\textup{Family} \(1\) has one orbit.

\item
\textup{Family} \(11\) has one orbit.

\item
\textup{Family} \(1'\) has one orbit when \(r'\ne2\).
When \(r'=2\), there are two orbits.
In this case,
\[
\mathbb R^4
=
U_{\mathrm{triv}}\oplus U_{\mathrm{nontriv}}.
\]
The group \(G\) acts trivially on \(U_{\mathrm{triv}}\) and by the nontrivial real character on \(U_{\mathrm{nontriv}}\).
A line in \(U_{\mathrm{triv}}\) and a line in \(U_{\mathrm{nontriv}}\) represent the two orbits.

\item
For \textup{Family} \(11'\) with \(r'>2\), write \(\ell'=(1-s)/r'\) when \(G_{1'}\) belongs to \familyref{1'-6} and \(\ell'=(1+s)/r'\) when \(G_{1'}\) belongs to \familyref{1'-7}.
There are two orbits when \(\ell'\) is even and one orbit when \(\ell'\) is odd.
If \(\ell'\) is even and \(G_{1'}\) belongs to \familyref{1'-6}, the two orbits are represented by \(\operatorname{span}\{1\}\) and \(\operatorname{span}\{\mathbf{i}\}\).
If \(\ell'\) is even and \(G_{1'}\) belongs to \familyref{1'-7}, the two orbits are represented by \(\operatorname{span}\{\mathbf{j}\}\) and \(\operatorname{span}\{\mathbf{k}\}\).
When \(r'=2\), there are two orbits.
The unique line fixed pointwise by \(G\) represents one orbit.
The remaining three invariant lines form the other orbit.
\end{enumerate}

Orthogonal complementation of \(G\)-invariant lines induces a bijection between these orbits and the \(N^+_\rho\)-orbits of invariant \(3\)-dimensional subspaces.
\end{proposition}

\begin{proof}
The proof of Lemma~\ref{lem:coordinate-factor-swap} shows that conjugation by \(\beta\) replaces \(s\) by \(-s\).
For \textup{Families} \(1\) and \(11\), conjugation by \(\beta\) induces a bijection between the normalizer orbits of invariant lines under the conditions defining \familyref{1-6} and those under the conditions defining \familyref{1-7}.
For \textup{Families} \(1'\) and \(11'\), conjugation by \(\beta\) induces a bijection between the normalizer orbits of invariant lines under the conditions defining \familyref{1'-6} and those under the conditions defining \familyref{1'-7}.
For each of \textup{Families} \(1\), \(11\), \(1'\), and \(11'\), the map \(\beta\) conjugates the representatives in the two parameter ranges but need not normalize either representative.
It suffices to consider \familyref{1-6} and \familyref{1'-6}.

Let \(U\) be a \(G\)-invariant \(2\)-dimensional subspace.
Assume that \(g|_U\in\SO(U)\) for every \(g\in G\).
Then \(R\oplus I_{U^\perp}\in C_{\SO(4)}(G)\) for every \(R\in\SO(U)\).
All lines in \(U\) lie in one \(N^+_\rho\)-orbit.

Suppose that \(G\) belongs to \familyref{1-6}.
Then \(m=n=1\), \(\gamma_1=\gamma_2=-I_4\), and
\[
\gamma_3
=
(-1)^{l_-}I_{\operatorname{span}\{1,\mathbf{i}\}}
\oplus
(-1)^{l_-}\operatorname{Rot}(2\pi/r).
\]
By Lemma~\ref{lem:rot-real-eigen}, every invariant line is contained in \(\operatorname{span}\{1,\mathbf{i}\}\) when \(r>2\).
Since \(g|_{\operatorname{span}\{1,\mathbf{i}\}}\in\SO(\operatorname{span}\{1,\mathbf{i}\})\) for every \(g\in G\), all invariant lines lie in one \(N^+_\rho\)-orbit.
For \(r=1\), the action on \(\mathbb R^4\) is scalar and \(C_{\SO(4)}(G)=\SO(4)\).
For \(r=2\), every invariant line lies in \(\operatorname{span}\{1,\mathbf{i}\}\) or \(\operatorname{span}\{\mathbf{j},\mathbf{k}\}\).
For each \(U\in\{\operatorname{span}\{1,\mathbf{i}\},\operatorname{span}\{\mathbf{j},\mathbf{k}\}\}\), all lines in \(U\) lie in a single \(N^+_\rho\)-orbit.
By Lemma~\ref{lem:coordinate-factor-swap}, an element of \(N^+_\rho\) interchanges these two subspaces.
Thus \textup{Family} \(1\) has one orbit.

Suppose that \(G\) belongs to \familyref{1'-6}.
Then \(m=n=1\), \(\gamma_1=\gamma_2=I_4\), and
\[
\gamma_3
=
(-1)^{l'_-}I_{\operatorname{span}\{1,\mathbf{i}\}}
\oplus
(-1)^{l'_-}\operatorname{Rot}(2\pi/r').
\]
By Lemma~\ref{lem:rot-real-eigen}, every invariant line is contained in \(\operatorname{span}\{1,\mathbf{i}\}\) when \(r'>2\).
Since \(g|_{\operatorname{span}\{1,\mathbf{i}\}}\in\SO(\operatorname{span}\{1,\mathbf{i}\})\) for every \(g\in G\), all invariant lines lie in one \(N^+_\rho\)-orbit.
For \(r'=1\), the condition \(\gcd(s,2r')=1\) implies that \(l'_-=1-s\) is even.
The representation is trivial, and \(\SO(4)\) acts transitively on the lines.
For \(r'=2\), write \(U_{\mathrm{triv}}\) and \(U_{\mathrm{nontriv}}\) for the two \(2\)-dimensional summands on which \(G\) acts by the trivial and nontrivial real characters, respectively.
Every invariant line is contained in \(U_{\mathrm{triv}}\) or \(U_{\mathrm{nontriv}}\).
The rotations in either summand centralize \(G\), so the lines in each summand form one \(N^+_\rho\)-orbit.
Every element of \(N_{\SO(4)}(G)\) preserves \(\operatorname{Fix}(G)=U_{\mathrm{triv}}\).
The two orbits are therefore distinct.

By Proposition~\ref{prop:criterion-line-summand} and Definition~\ref{def:representative-family-parameters}, \(m=n=1\) and \(r\ge2\) in \textup{Family} \(11\).
By Proposition~\ref{prop:case-family-11},
\[
G_{11}=\langle G_1,\gamma_4\rangle,
\qquad
\gamma_4=F\oplus F,
\qquad
F=
\begin{bmatrix}
1&0\\
0&-1
\end{bmatrix}.
\]
For \(r>2\), the invariant lines are \(\operatorname{span}\{1\}\) and \(\operatorname{span}\{\mathbf{i}\}\).
Set \(J=\operatorname{Rot}(\pi/2)\) and \(h=J\oplus J\).
The map \(h\) commutes with \(G_1\).
Since \(JFJ^{-1}=-F\), we have \(h\gamma_4h^{-1}=(-I_4)\gamma_4=\gamma_1\gamma_4\).
Thus \(h\in N^+_\rho\), and \(h\) interchanges the two lines.

For \(r=2\),
\[
G_{11}
=
\left\{
\operatorname{diag}(\varepsilon_1,\varepsilon_2,\varepsilon_3,\varepsilon_4)
\mathrel{}\middle|\mathrel{}
\varepsilon_i\in\{\pm1\}\text{ for }1\le i\le4,
\quad
\varepsilon_1\varepsilon_2\varepsilon_3\varepsilon_4=1
\right\}.
\]
The lines \(\operatorname{span}\{1\}\), \(\operatorname{span}\{\mathbf{i}\}\), \(\operatorname{span}\{\mathbf{j}\}\), and \(\operatorname{span}\{\mathbf{k}\}\) are the one-dimensional character spaces of \(G_{11}\).
Signed permutations of the coordinates that preserve orientation normalize \(G_{11}\) and act transitively on these lines.
Thus \textup{Family} \(11\) has one orbit.

By Proposition~\ref{prop:criterion-line-summand} and Definition~\ref{def:representative-family-parameters}, \(m=n=1\) and \(r'\ge2\) in \textup{Family} \(11'\).
By Proposition~\ref{prop:case-family-11p},
\[
G_{11'}=\langle G_{1'},\gamma_4\rangle,
\qquad
\gamma_4=F\oplus F.
\]

Assume that \(r'>2\).
If \(G_{1'}\) belongs to \familyref{1'-6}, retain the chosen representative.
If \(G_{1'}\) belongs to \familyref{1'-7}, conjugate the chosen representative by \(\beta\).
This replaces \(s\) by \(-s\).
After this conjugation, write \(s\) for the resulting parameter.
Then \(\ell'=(1-s)/r'\).
The invariant lines and the third generator are
\[
L_0=\operatorname{span}\{1\},
\qquad
L_1=\operatorname{span}\{\mathbf{i}\},
\qquad
\gamma_3
=
(-1)^{\ell'}I_{\operatorname{span}\{1,\mathbf{i}\}}
\oplus
(-1)^{\ell'}\operatorname{Rot}(2\pi/r').
\]

Suppose that \(\ell'\) is even.
All generators fix \(L_0\) pointwise.
The generator \(\gamma_4\) acts nontrivially on \(L_1\).
By Lemma~\ref{lem:rot-real-eigen}, \(\operatorname{span}\{\mathbf{j},\mathbf{k}\}\) contains no nonzero fixed vector.
Thus \(\operatorname{Fix}(G_{11'})=L_0\).
Since the normalizer preserves \(L_0\), the lines \(L_0\) and \(L_1\) lie in distinct orbits.

Suppose that \(\ell'\) is odd.
Set
\[
h
=
\operatorname{Rot}(\pi/2)
\oplus
\operatorname{Rot}(\pi/2+\pi/r').
\]
Then \(h\in\SO(4)\), and \(h\) commutes with \(\gamma_3\).
For \(\theta\in\mathbb R\),
\[
\operatorname{Rot}(\theta)
F
\operatorname{Rot}(-\theta)
=
\operatorname{Rot}(2\theta)F.
\]
It follows that \(h\gamma_4h^{-1}=\gamma_3\gamma_4\).
Thus \(h\in N^+_\rho\), and \(h\) interchanges \(L_0\) and \(L_1\).

For parameters satisfying \familyref{1'-7}, the two invariant lines in the original coordinates are \(\operatorname{span}\{\mathbf{j}\}\) and \(\operatorname{span}\{\mathbf{k}\}\).

For \(r'=2\), after an orientation-preserving change of coordinates,
\[
G_{11'}
=
\left\{
\operatorname{diag}(1,\varepsilon_2,\varepsilon_3,\varepsilon_4)
\mathrel{}\middle|\mathrel{}
\varepsilon_i\in\{\pm1\}\text{ for }2\le i\le4,
\quad
\varepsilon_2\varepsilon_3\varepsilon_4=1
\right\}.
\]
The four coordinate axes are the distinct one-dimensional character spaces of \(G_{11'}\).
The first axis is the unique line fixed pointwise.
Signed permutations of the coordinates that preserve orientation and fix the first axis normalize \(G_{11'}\) and act transitively on the other three axes.
The fixed line forms one \(N^+_\rho\)-orbit, and the remaining three character spaces form the other.

For every \(T\in N^+_\rho\) and every invariant line \(L\), orthogonality gives \(T(L^\perp)=T(L)^\perp\).
Orthogonal complementation therefore induces a bijection between the orbits of invariant lines and the orbits of invariant \(3\)-dimensional subspaces.
\end{proof}

\subsection{\textup{Families} \texorpdfstring{\(11\), \(11'\), \(3\), and \(34\)}{11, 11', 3, and 34}}\label{subsec:remaining-plane-orbits}

\begin{proposition}\label{prop:pillow-index-two-handle-orbits}
Let \(G<\SO(4)\) be the chosen representative in \textup{Family} \(11\) or \(11'\).
Whenever they occur, \(m,n,r,s,r'\) are the parameters of \(G\) in Definition~\ref{def:representative-family-parameters}.
Let \(\rho\colon G\hookrightarrow\SO(4)\) be the inclusion.
Every invariant \(2\)-dimensional subspace lies in an \(N^+_\rho\)-orbit represented by \(\operatorname{span}\{1,\mathbf{i}\}\) or \(\operatorname{span}\{\mathbf{j},\mathbf{k}\}\).
The orbits are classified as follows.

\begin{enumerate}[label=\textup{(\roman*)}]
\item
For \textup{Family} \(11\), the two representatives lie in the same orbit precisely when \(r=1\) or \(r=2\).
For \(r>2\), they represent two distinct orbits.

\item
For \textup{Family} \(11'\), the two representatives lie in the same orbit precisely when \(r'=1\).
For \(r'>1\), they represent two distinct orbits.
When \(m=n=1\) and \(r'=2\), the two orbits are distinguished by whether the \(2\)-dimensional subspace contains the unique line fixed pointwise by \(G\).
\end{enumerate}

Orthogonal complementation interchanges the two orbits whenever they are distinct.
\end{proposition}

\begin{proof}
The restrictions \(mr,nr\ge2\) exclude the conditions defining \familyref{1-1} and \familyref{1-3}.
The restrictions \(mr',nr'\ge2\) exclude the conditions defining \familyref{1'-1} and \familyref{1'-3}.
The calculation in the proof of Proposition~\ref{prop:cyclic-index-two-handle-orbits} shows that there is no nonzero equivariant linear map between \(\operatorname{span}\{1,\mathbf{i}\}\) and \(\operatorname{span}\{\mathbf{j},\mathbf{k}\}\).

Except when \((m,n,r)=(1,1,2)\), these are the only invariant \(2\)-dimensional subspaces for \(G_1\).
Since \(G_1\subset G_{11}\) and \(\gamma_4\) preserves both subspaces, they are also the only invariant \(2\)-dimensional subspaces for \(G_{11}\) outside this case.

Except when \((m,n,r')=(1,1,2)\), these are the only invariant \(2\)-dimensional subspaces for \(G_{1'}\).
Since \(G_{1'}\subset G_{11'}\) and \(\gamma_4\) preserves both subspaces, they are also the only invariant \(2\)-dimensional subspaces for \(G_{11'}\) outside this case.

By Lemma~\ref{lem:coordinate-factor-swap}, \(\operatorname{span}\{1,\mathbf{i}\}\) and \(\operatorname{span}\{\mathbf{j},\mathbf{k}\}\) lie in the same \(N^+_\rho\)-orbit precisely when \(r\in\{1,2\}\) in \textup{Family} \(11\) and when \(r'=1\) in \textup{Family} \(11'\).
Orthogonal complementation interchanges these two subspaces.

If \(m=n=1\) and \(r=2\), the lines \(\operatorname{span}\{1\}\), \(\operatorname{span}\{\mathbf{i}\}\), \(\operatorname{span}\{\mathbf{j}\}\), and \(\operatorname{span}\{\mathbf{k}\}\) are the distinct one-dimensional character spaces of \(G_{11}\).
Every invariant \(2\)-dimensional subspace is a sum of two of these lines.
Signed permutations of the coordinates that preserve orientation normalize \(G_{11}\) and act transitively on the two-element subsets of these four lines.
Thus all invariant \(2\)-dimensional subspaces lie in one \(N^+_\rho\)-orbit.

If \(m=n=1\) and \(r'=2\), the lines \(\operatorname{span}\{1\}\), \(\operatorname{span}\{\mathbf{i}\}\), \(\operatorname{span}\{\mathbf{j}\}\), and \(\operatorname{span}\{\mathbf{k}\}\) are the distinct character spaces of \(G_{11'}\).
Exactly one of these lines is fixed pointwise.
Every invariant \(2\)-dimensional subspace is a sum of two of these lines.
The normalizer preserves the fixed line.
Signed permutations of the coordinates that preserve orientation and fix this line act transitively on the other three lines.
There are exactly two orbits, distinguished by whether the subspace contains the fixed line.
Orthogonal complementation interchanges these two orbits.
\end{proof}

\begin{proposition}\label{prop:block-swap-two-plane-handle-orbits}
Let \(G<\SO(4)\) be the chosen representative in \textup{Family} \(3\) with \(m=1\) and \(n\ge2\), or in \textup{Family} \(34\) with \(m=1\) and odd \(n\ge3\).
Here \(m,n\) are the parameters of \(G\) in Definition~\ref{def:representative-family-parameters}.
Let \(\rho\colon G\hookrightarrow\SO(4)\) be the inclusion.
Then \(\{W_\theta\mid\theta\in\mathbb R\}\) is a single \(N^+_\rho\)-orbit.
Orthogonal complementation leaves this orbit invariant.
\end{proposition}

\begin{proof}
For \(t\in\mathbb R\), set \(T_t=\operatorname{diag}(\operatorname{Rot}(t),\operatorname{Rot}(t))\).
By the formulas for the generators in the proofs of Propositions~\ref{prop:case-family-3} and~\ref{prop:case-family-34}, \(T_t\in C_{\SO(4)}(G)\) when \(m=1\).
For \(\theta,t\in\mathbb R\), the action satisfies \(T_t(W_\theta)=W_{\theta+2t}\).
Thus the subspaces \(W_\theta\), \(\theta\in\mathbb R\), form a single \(N^+_\rho\)-orbit.

By Lemma~\ref{lem:special-block-switch}, \(W_\theta^\perp=W_{\theta+\pi}\) for every \(\theta\in\mathbb R\).
Orthogonal complementation leaves this orbit invariant.
\end{proof}

\begin{proof}[\textup{Proof of Theorem B}]
By Theorem~\ref{thm:4d-morse-index-classification}, the possible dimensions of invariant subspaces are those in Table~\ref{tab:normalizer-orbit-representatives}.
For the trivial group, \(N_{\SO(4)}(G)=\SO(4)\) acts transitively on \(\operatorname{Gr}_k(\mathbb R^4)\) for every \(k\).
The only subspaces of dimensions \(0\) and \(4\) are \(0\) and \(\mathbb H\), respectively.

Proposition~\ref{prop:line-handle-orbits} determines the orbits of invariant lines in \textup{Families} \(1\), \(1'\), \(11\), and \(11'\).
For \textup{Family} \(1\), the divisibility conditions are \(r\mid(1-s)\) and \(r\mid(1+s)\).
For \textup{Family} \(1'\), the divisibility conditions are \(r'\mid(1-s)\) and \(r'\mid(1+s)\).
The proof of Proposition~\ref{prop:line-handle-orbits} identifies the representatives in the chosen coordinates.
For \textup{Family} \(11'\) with \(m=n=1\) and \(r'=2\), the unique line fixed pointwise is \(\operatorname{span}\{1\}\) when \((1-s)/2\) is even.
It is \(\operatorname{span}\{\mathbf{j}\}\) when \((1-s)/2\) is odd.
In either case, \(\operatorname{span}\{1\}\) and \(\operatorname{span}\{\mathbf{j}\}\) represent the two orbits.
Corollary~\ref{cor:unique-index-one-three-TOI} determines the orbits of invariant lines in \textup{Families} \(21\), \(21'\), \(26\), \(26'\), \(26''\), \(31\), and \(31'\).
Orthogonal complementation is \(N^+_\rho\)-equivariant and determines the representatives for index \(3\).

Propositions~\ref{prop:cyclic-index-two-handle-orbits} and~\ref{prop:pillow-index-two-handle-orbits} determine the representatives for index \(2\) in \textup{Families} \(1\), \(1'\), \(11\), and \(11'\).
For \textup{Family} \(1\) with \(m=n=1\) and \(r=2\), and for \textup{Family} \(1'\) with \(m=n=1\) and \(r'=2\), choose \(\ell_1=\operatorname{span}\{1\}\) and \(\ell_2=\operatorname{span}\{\mathbf{j}\}\).
Then \(\ell_1\oplus\ell_2=\operatorname{span}\{1,\mathbf{j}\}\).
By Proposition~\ref{prop:block-swap-two-plane-handle-orbits}, \(W_0\) represents the unique orbit in \textup{Families} \(3\) and \(34\) with \(m=1\).
By Definition~\ref{def:block-graph-plane}, \(W_0=\operatorname{span}\{1+\mathbf{j},\mathbf{i}-\mathbf{k}\}\).
Under factor interchange, the representatives of the orbits are obtained by applying \(V\mapsto\mu V\).

By Proposition~\ref{prop:ordered-handle-orbits}, these orbits are in bijection with the isomorphism classes of ordered handles.
Lemma~\ref{lem:orbifold-morse-normal-form} identifies the invariant splitting at each nondegenerate critical point on an orbifold.
Definition~\ref{def:ordered-handle} defines the ordered handle determined by this splitting.
For each representative \(V\), set \(V^+=V^\perp\).
By Lemma~\ref{lem:main-hessian-splitting-principle}, the invariant quadratic function \(Q_{V,V^+}\) realizes the splitting \(V\oplus V^+\).
\end{proof}

\appendix

\section{Calculations of invariant subspaces from generators}\label{app:morse-index-case-analysis}

In this appendix, we compute the invariant subspaces directly from generators for the families not covered by Proposition~\ref{prop:criterion-q8}.
We use the representatives in the completed Du Val list in Section~\ref{subsec:so4-finite-subgroups}~\cite{val1964homographies,mecchia2015fibered}.
The proofs use Proposition~\ref{prop:generator} and the linear algebra lemmas in Section~\ref{subsec:q8-criterion}.
The subfamily labels are those of Definition~\ref{def:parameter-cases-1-1p}.
We use the notation \(M_\theta\) and \(W_\theta\) from Definition~\ref{def:block-graph-plane}.
Sections~\ref{subsec:morse-index-theorem} and~\ref{sec:four-dimensional-handle-classification} use these calculations.
The proofs do not use the classification of Morse indices in Section~\ref{subsec:morse-index-theorem} or the classification of ordered handles in Section~\ref{sec:four-dimensional-handle-classification}.

Section~\ref{subsec:appendix-general-lemmas} proves the lemmas on invariant subspaces used in the calculations.
Section~\ref{subsec:appendix-cyclic-dihedral} determines the invariant subspaces in \textup{Families} \(1\), \(1'\), \(3\), \(11\), and \(11'\), including the reducible case \(m=1\) of \textup{Family} \(3\).
Section~\ref{subsec:appendix-polyhedral} determines the invariant lines and irreducible \(3\)-dimensional complements in \textup{Families} \(21\), \(21'\), \(26\), \(26'\), \(26''\), \(31\), and \(31'\).
It also proves irreducibility for \textup{Families} \(32\) and \(32'\).
Section~\ref{subsec:appendix-thirty-three-four} determines the invariant subspaces in \textup{Families} \(33\), \(33'\), and \(34\).
It also corrects the quotient isomorphism for \textup{Family} \(33'\).
Under factor interchange, \(V\mapsto\mu V\) maps the invariant subspaces and normalizer orbits of each chosen representative to those of the group obtained by interchanging the factors.

Throughout the appendix, write
\[
\mathbb H
=
\mathbb C\oplus\mathbb C\mathbf{j}
=
\operatorname{span}\{1,\mathbf{i}\}
\oplus
\operatorname{span}\{\mathbf{j},\mathbf{k}\}.
\]

\subsection{General lemmas on invariant subspaces}\label{subsec:appendix-general-lemmas}

\begin{lemma}\label{lem:cyclic-block-line-test}
Let \(a,b,c_-,c_+\in\mathbb R\).
Let \(G<\SO(4)\) be generated by
\[
\begin{aligned}
\gamma_1&=\operatorname{diag}(\operatorname{Rot}(a),\operatorname{Rot}(a)),\\
\gamma_2&=\operatorname{diag}(\operatorname{Rot}(-b),\operatorname{Rot}(b)),\\
\gamma_3&=\operatorname{diag}(\operatorname{Rot}(c_-),\operatorname{Rot}(c_+)).
\end{aligned}
\]
If either \(a\not\equiv0,\pi\pmod{2\pi}\) or \(b\not\equiv0,\pi\pmod{2\pi}\), then \(G\) has no invariant line.

Suppose that \(\gamma_1\) and \(\gamma_2\) are scalar.
Then \(G\) has an invariant line precisely when
\[
c_-\equiv0,\pi\pmod{2\pi}
\qquad\text{or}\qquad
c_+\equiv0,\pi\pmod{2\pi}.
\]
If \(c_-\equiv0,\pi\pmod{2\pi}\), then \(\operatorname{span}\{1,\mathbf{i}\}\) is a direct sum of invariant lines.
If \(c_+\equiv0,\pi\pmod{2\pi}\), then \(\operatorname{span}\{\mathbf{j},\mathbf{k}\}\) is a direct sum of invariant lines.
\end{lemma}

\begin{proof}
By Lemma~\ref{lem:rot-real-eigen}, \(\gamma_1\) has no real eigenvector if \(a\not\equiv0,\pi\pmod{2\pi}\).
The transformation \(\gamma_2\) has no real eigenvector if \(b\not\equiv0,\pi\pmod{2\pi}\).
Under either condition, \(G\) has no invariant line.

Suppose that \(\gamma_1\) and \(\gamma_2\) are scalar.
Then a line is \(G\)-invariant precisely when it is invariant under \(\gamma_3\).
By Lemma~\ref{lem:rot-real-eigen}, the transformation \(\gamma_3\) has a real eigenvector precisely when \(c_-\equiv0,\pi\pmod{2\pi}\) or \(c_+\equiv0,\pi\pmod{2\pi}\).
If \(c_-\equiv0,\pi\pmod{2\pi}\), then every line in \(\operatorname{span}\{1,\mathbf{i}\}\) is invariant.
If \(c_+\equiv0,\pi\pmod{2\pi}\), then every line in \(\operatorname{span}\{\mathbf{j},\mathbf{k}\}\) is invariant.
\end{proof}

\begin{lemma}\label{lem:two-block-swap-irreducible}
Let \(\alpha,\delta\in\mathbb R\) satisfy \(\alpha\not\equiv0,\pi\pmod{2\pi}\) and \(\delta\not\equiv0,\pi\pmod{2\pi}\).
Let \(G<\SO(4)\) contain
\[
A=\operatorname{diag}(\operatorname{Rot}(\alpha),\operatorname{Rot}(\alpha)),
\qquad
B=\operatorname{diag}(\operatorname{Rot}(-\delta),\operatorname{Rot}(\delta)).
\]
Assume that \(G\) contains an element \(S\) satisfying \(S(\operatorname{span}\{1,\mathbf{i}\})=\operatorname{span}\{\mathbf{j},\mathbf{k}\}\).
Then the real \(G\)-module \(\mathbb R^4\) is irreducible.
\end{lemma}

\begin{proof}
By Lemma~\ref{lem:rot-real-eigen}, \(A\) has no real eigenvector.
Thus \(G\) has no invariant line.
The orthogonal complement of a \(G\)-invariant \(3\)-dimensional subspace is a \(G\)-invariant line.
Thus \(G\) has no invariant \(3\)-dimensional subspace.

Let \(P\) be a \(G\)-invariant \(2\)-dimensional subspace.
By Lemma~\ref{lem:matrix-commute} applied to \(A\), either \(P=\operatorname{span}\{\mathbf{j},\mathbf{k}\}\) or \(P=\operatorname{graph}(T)\) for a linear map \(T\colon \operatorname{span}\{1,\mathbf{i}\}\to \operatorname{span}\{\mathbf{j},\mathbf{k}\}\).
Suppose that \(P=\operatorname{graph}(T)\).
By Lemma~\ref{lem:rot-invar-basis}, invariance of \(P\) under \(A\) implies that \(T\) is complex-linear.
Invariance of \(P\) under \(B\) implies that \(T\) is conjugate-linear.
It follows that \(T=0\), so \(P=\operatorname{span}\{1,\mathbf{i}\}\).

The only \(2\)-dimensional subspaces invariant under both \(A\) and \(B\) are \(\operatorname{span}\{1,\mathbf{i}\}\) and \(\operatorname{span}\{\mathbf{j},\mathbf{k}\}\).
Since \(S\) interchanges these two subspaces, neither is \(G\)-invariant.
Therefore the real \(G\)-module \(\mathbb R^4\) is irreducible.
\end{proof}

Set
\[
J_0=
\begin{bmatrix}
0&-1\\
1&0
\end{bmatrix}.
\]

\begin{lemma}\label{lem:special-block-switch}
Let \(\delta\in\mathbb R\) with \(\delta\not\equiv0,\pi\pmod{2\pi}\).
Set \(B=\operatorname{diag}(\operatorname{Rot}(-\delta),\operatorname{Rot}(\delta))\).
Define \(S\colon\mathbb R^2\oplus\mathbb R^2\to\mathbb R^2\oplus\mathbb R^2\) by \(S(v,w)=(J_0w,-J_0v)\).
Set \(G=\langle B,S\rangle\).
Then \(G\) has no invariant line or invariant \(3\)-dimensional subspace.
Its invariant \(2\)-dimensional subspaces are precisely \(W_\theta\) for \(\theta\in\mathbb R\).
For every \(\theta\in\mathbb R\), \(W_\theta^\perp=W_{\theta+\pi}\).
\end{lemma}

\begin{proof}
By Lemma~\ref{lem:rot-real-eigen}, \(B\) has no real eigenvector.
Thus \(G\) has no invariant line.
The orthogonal complement of a \(G\)-invariant \(3\)-dimensional subspace is a \(G\)-invariant line.
Thus \(G\) has no invariant \(3\)-dimensional subspace.

Let \(P\) be a \(G\)-invariant \(2\)-dimensional subspace.
Since \(S\) interchanges \(\operatorname{span}\{1,\mathbf{i}\}\) and \(\operatorname{span}\{\mathbf{j},\mathbf{k}\}\), neither subspace is \(G\)-invariant.
By Lemma~\ref{lem:matrix-commute} applied to \(B\), \(P=\operatorname{graph}(T)\) for a linear map \(T\colon \operatorname{span}\{1,\mathbf{i}\}\to \operatorname{span}\{\mathbf{j},\mathbf{k}\}\).
By Lemma~\ref{lem:rot-invar-basis}, there exist \(p,q\in\mathbb R\) such that
\[
T=
\begin{bmatrix}
p&q\\
q&-p
\end{bmatrix}.
\]
Invariance of \(P\) under \(S\) is equivalent to \(TJ_0T=-J_0\).
Since \(TJ_0T=-(p^2+q^2)J_0\), this condition is equivalent to \(p^2+q^2=1\).
Thus \(P=W_\theta\) for some \(\theta\in\mathbb R\).

Since \(M_\theta^*=M_\theta\) and \(M_\theta^2=I_2\), \(W_\theta^\perp=\{(v,-M_\theta v)\mid v\in\mathbb R^2\}=W_{\theta+\pi}\).
\end{proof}

Set \(w_1=\frac12(1+\mathbf{i}+\mathbf{j}+\mathbf{k})\) and \(W=\operatorname{Im}\mathbb H=\operatorname{span}\{\mathbf{i},\mathbf{j},\mathbf{k}\}\).
Set \(H_T=\Phi(\{(q,q)\mid q\in T^*\})\).
Then \(H_T\cong T^*/\{\pm1\}\cong\mathbf T\).

\begin{lemma}\label{lem:diagonal-polyhedral-module}
Let \(P\in\{T^*,O^*,I^*\}\).
Let \(\widetilde G<S^3\times S^3\) be a finite subgroup satisfying
\[
\{(q,q)\mid q\in T^*\}
\subset \widetilde G
\subset
\{(p,\varepsilon p)\mid p\in P,\ \varepsilon\in\{\pm1\}\}.
\]
Set \(G=\Phi(\widetilde G)\).
Then the splitting
\[
\mathbb H=\operatorname{span}\{1\}\oplus W
\]
is \(G\)-invariant.
The real \(G\)-module \(W\) is irreducible.
\end{lemma}

\begin{proof}
For each \(g\in G\), choose \(p\in P\) and \(\varepsilon\in\{\pm1\}\) such that \(g=\Phi(p,\varepsilon p)\).
Then \(g(v)=\varepsilon pvp^{-1}\).
Thus the splitting \(\mathbb H=\operatorname{span}\{1\}\oplus W\) is \(G\)-invariant.

For \(q\in T^*\), the element \(\Phi(q,q)\) acts on \(W\) by \(v\mapsto qvq^{-1}\).
Conjugation by \(w_1\) has the unique invariant line \(\operatorname{span}\{\mathbf{i}+\mathbf{j}+\mathbf{k}\}\) in \(W\).
Conjugation by \(\mathbf{i}\) does not preserve this line.
Therefore \(H_T\) has no invariant line in \(W\).

If \(H_T\) had an invariant \(2\)-dimensional subspace of \(W\), its orthogonal complement in \(W\) would be an \(H_T\)-invariant line.
Hence \(W\) is an irreducible real \(H_T\)-module.
Since \(H_T\subset G\), the real \(G\)-module \(W\) is irreducible.
\end{proof}

\subsection{\textup{Families} \texorpdfstring{\(1\), \(1'\), \(3\), \(11\), and \(11'\)}{1, 1', 3, 11, and 11'}}\label{subsec:appendix-cyclic-dihedral}

For \textup{Family} \(1\), set
\[
\widetilde G_1
=
(\Cyc{2mr}/\Cyc{2m},\Cyc{2nr}/\Cyc{2n})_s,
\qquad
G_1=\Phi(\widetilde G_1).
\]

\begin{proposition}\label{prop:case-family-1}
Let \(m,n,r,s\) be the family parameters of \(G_1\) as in Definition~\ref{def:representative-family-parameters}.
The \(2\)-dimensional subspaces \(\operatorname{span}\{1,\mathbf{i}\}\) and \(\operatorname{span}\{\mathbf{j},\mathbf{k}\}\) are \(G_1\)-invariant.
The following hold.

\begin{enumerate}[label=\textup{(\roman*)}]
\item
If \(G_1\) belongs to neither \familyref{1-6} nor \familyref{1-7}, both subspaces are irreducible real \(G_1\)-modules.
They form a splitting of type \(2\oplus2\).

\item
If \(G_1\) belongs to exactly one of \familyref{1-6} and \familyref{1-7}, then \(\mathbb R^4\) admits an invariant splitting of type \(1\oplus1\oplus2\).

\item
If \(G_1\) belongs to both \familyref{1-6} and \familyref{1-7}, then \(\mathbb R^4\) admits an invariant splitting of type \(1\oplus1\oplus1\oplus1\).
\end{enumerate}
\end{proposition}

\begin{proof}
The generators of \(G_1\) in Section~\ref{subsec:parameter-cases} preserve \(\operatorname{span}\{1,\mathbf{i}\}\) and \(\operatorname{span}\{\mathbf{j},\mathbf{k}\}\).

By Lemma~\ref{lem:cyclic-block-line-test}, \(G_1\) has an invariant line only if \(m=n=1\).
In this case, the restriction of \(\gamma_3\) to \(\operatorname{span}\{1,\mathbf{i}\}\) is scalar precisely when \((1-s)/r\in\mathbb Z\).
The restriction of \(\gamma_3\) to \(\operatorname{span}\{\mathbf{j},\mathbf{k}\}\) is scalar precisely when \((1+s)/r\in\mathbb Z\).
These conditions define \familyref{1-6} and \familyref{1-7}, respectively.

If neither condition holds, neither \(2\)-dimensional subspace contains an invariant line, so both subspaces are irreducible.
If exactly one condition holds, the subspace on which \(\gamma_3\) is scalar splits into invariant lines.
If both conditions hold, both \(2\)-dimensional subspaces split into invariant lines.
\end{proof}

For \textup{Family} \(1'\), write \(r=2r'\).
Set
\[
\widetilde G_{1'}
=
(\Cyc{mr}/\Cyc m,\Cyc{nr}/\Cyc n)_s,
\qquad
G_{1'}=\Phi(\widetilde G_{1'}).
\]

\begin{proposition}\label{prop:case-family-1p}
Let \(m,n,r',s\) be the family parameters of \(G_{1'}\) as in Definition~\ref{def:representative-family-parameters}.
The \(2\)-dimensional subspaces \(\operatorname{span}\{1,\mathbf{i}\}\) and \(\operatorname{span}\{\mathbf{j},\mathbf{k}\}\) are \(G_{1'}\)-invariant.
The following hold.

\begin{enumerate}[label=\textup{(\roman*)}]
\item
If \(G_{1'}\) belongs to neither \familyref{1'-6} nor \familyref{1'-7}, both subspaces are irreducible real \(G_{1'}\)-modules.
They form a splitting of type \(2\oplus2\).

\item
If \(G_{1'}\) belongs to exactly one of \familyref{1'-6} and \familyref{1'-7}, then \(\mathbb R^4\) admits an invariant splitting of type \(1\oplus1\oplus2\).

\item
If \(G_{1'}\) belongs to both \familyref{1'-6} and \familyref{1'-7}, then \(\mathbb R^4\) admits an invariant splitting of type \(1\oplus1\oplus1\oplus1\).
\end{enumerate}
\end{proposition}

\begin{proof}
The generators of \(G_{1'}\) in Section~\ref{subsec:parameter-cases} preserve \(\operatorname{span}\{1,\mathbf{i}\}\) and \(\operatorname{span}\{\mathbf{j},\mathbf{k}\}\).

By Definition~\ref{def:representative-family-parameters}, \(m\) and \(n\) are odd.
By Lemma~\ref{lem:cyclic-block-line-test}, \(G_{1'}\) has an invariant line only if \(m=n=1\).
In this case, the restriction of \(\gamma_3\) to \(\operatorname{span}\{1,\mathbf{i}\}\) is scalar precisely when \((1-s)/r'\in\mathbb Z\).
The restriction of \(\gamma_3\) to \(\operatorname{span}\{\mathbf{j},\mathbf{k}\}\) is scalar precisely when \((1+s)/r'\in\mathbb Z\).
These conditions define \familyref{1'-6} and \familyref{1'-7}, respectively.

If neither condition holds, neither \(2\)-dimensional subspace contains an invariant line, so both subspaces are irreducible.
If exactly one condition holds, the subspace on which \(\gamma_3\) is scalar splits into invariant lines.
If both conditions hold, both \(2\)-dimensional subspaces split into invariant lines.
\end{proof}

For \textup{Family} \(11\), set
\[
\widetilde G_{11}
=
(\BDih{4mr}/\Cyc{2m},\BDih{4nr}/\Cyc{2n})_s,
\qquad
G_{11}=\Phi(\widetilde G_{11}).
\]

\begin{proposition}\label{prop:case-family-11}
Let \(m,n,r,s\) be the family parameters of \(G_{11}\) as in Definition~\ref{def:representative-family-parameters}.
Let \(G_1\) be the representative in \textup{Family} \(1\) with the same family parameters.
Set \(\gamma_4=\operatorname{diag}(1,-1,1,-1)\).
The group \(G_{11}\) is generated by \(G_1\) and \(\gamma_4\).
The following hold.

\begin{enumerate}[label=\textup{(\roman*)}]
\item
If \(G_1\) belongs to neither \familyref{1-6} nor \familyref{1-7}, then \(\mathbb R^4\) admits an invariant splitting of type \(2\oplus2\).

\item
If \(G_1\) belongs to exactly one of \familyref{1-6} and \familyref{1-7}, then \(\mathbb R^4\) admits an invariant splitting of type \(1\oplus1\oplus2\).

\item
If \(G_1\) belongs to both \familyref{1-6} and \familyref{1-7}, then \(\mathbb R^4\) admits an invariant splitting of type \(1\oplus1\oplus1\oplus1\).
\end{enumerate}
\end{proposition}

\begin{proof}
The element \((\mathbf{j},\mathbf{j})\) belongs to \(\widetilde G_{11}\).
We have \(\Phi(\mathbf{j},\mathbf{j})=\gamma_4\).
Thus \(G_{11}=\langle G_1,\gamma_4\rangle\).

The transformation \(\gamma_4\) preserves \(\operatorname{span}\{1,\mathbf{i}\}\) and \(\operatorname{span}\{\mathbf{j},\mathbf{k}\}\).
It preserves each of the lines \(\operatorname{span}\{1\}\), \(\operatorname{span}\{\mathbf{i}\}\), \(\operatorname{span}\{\mathbf{j}\}\), and \(\operatorname{span}\{\mathbf{k}\}\).
In the proof of Proposition~\ref{prop:case-family-1}, choose the coordinate axes whenever \(G_1\) acts by scalars on either \(2\)-dimensional coordinate subspace.
These choices preserve the dimensions of the summands.
Each splitting is then \(G_{11}\)-invariant.
\end{proof}

For \textup{Family} \(11'\), write \(r=2r'\).
Set
\[
\widetilde G_{11'}
=
(\BDih{2mr}/\Cyc m,\BDih{2nr}/\Cyc n)_s,
\qquad
G_{11'}=\Phi(\widetilde G_{11'}).
\]

\begin{proposition}\label{prop:case-family-11p}
Let \(m,n,r',s\) be the family parameters of \(G_{11'}\) as in Definition~\ref{def:representative-family-parameters}.
Let \(G_{1'}\) be the representative in \textup{Family} \(1'\) with the same family parameters.
Set \(\gamma_4=\operatorname{diag}(1,-1,1,-1)\).
The group \(G_{11'}\) is generated by \(G_{1'}\) and \(\gamma_4\).
The following hold.

\begin{enumerate}[label=\textup{(\roman*)}]
\item
If \(G_{1'}\) belongs to neither \familyref{1'-6} nor \familyref{1'-7}, then \(\mathbb R^4\) admits an invariant splitting of type \(2\oplus2\).

\item
If \(G_{1'}\) belongs to exactly one of \familyref{1'-6} and \familyref{1'-7}, then \(\mathbb R^4\) admits an invariant splitting of type \(1\oplus1\oplus2\).

\item
If \(G_{1'}\) belongs to both \familyref{1'-6} and \familyref{1'-7}, then \(\mathbb R^4\) admits an invariant splitting of type \(1\oplus1\oplus1\oplus1\).
\end{enumerate}
\end{proposition}

\begin{proof}
The element \((\mathbf{j},\mathbf{j})\) belongs to \(\widetilde G_{11'}\).
We have \(\Phi(\mathbf{j},\mathbf{j})=\gamma_4\).
Thus \(G_{11'}=\langle G_{1'},\gamma_4\rangle\).

The transformation \(\gamma_4\) preserves \(\operatorname{span}\{1,\mathbf{i}\}\) and \(\operatorname{span}\{\mathbf{j},\mathbf{k}\}\).
It preserves each of the lines \(\operatorname{span}\{1\}\), \(\operatorname{span}\{\mathbf{i}\}\), \(\operatorname{span}\{\mathbf{j}\}\), and \(\operatorname{span}\{\mathbf{k}\}\).
In the proof of Proposition~\ref{prop:case-family-1p}, choose the coordinate axes whenever \(G_{1'}\) acts by scalars on either \(2\)-dimensional coordinate subspace.
These choices preserve the dimensions of the summands.
Each splitting is then \(G_{11'}\)-invariant.
\end{proof}

For \textup{Family} \(3\), set
\[
\widetilde G_3
=
(\Cyc{4m}/\Cyc{2m},\BDih{4n}/\Cyc{2n}),
\qquad
G_3=\Phi(\widetilde G_3).
\]

\begin{proposition}\label{prop:case-family-3}
Let \(m,n\) be the family parameters of \(G_3\) as in Definition~\ref{def:representative-family-parameters}.
The following hold.

\begin{enumerate}[label=\textup{(\roman*)}]
\item
If \(m,n\ge2\), then the real \(G_3\)-module \(\mathbb H\) is irreducible.

\item
If \(m=1\), then \(G_3\) has no invariant line or invariant \(3\)-dimensional subspace.
In this case, its invariant \(2\)-dimensional subspaces are precisely \(W_\theta\) for \(\theta\in\mathbb R\).
\end{enumerate}
\end{proposition}

\begin{proof}
By Proposition~\ref{prop:generator}, \(\widetilde G_3\) is generated by
\[
(e^{\pi \mathbf{i}/m},1),
\qquad
(1,e^{\pi \mathbf{i}/n}),
\qquad
(e^{\pi \mathbf{i}/2m},\mathbf{j}).
\]
If \(m,n\ge2\), the images of these generators under \(\Phi\) satisfy the hypotheses of Lemma~\ref{lem:two-block-swap-irreducible} with \(\alpha=\pi/m\) and \(\delta=\pi/n\).
Hence the real \(G_3\)-module \(\mathbb H\) is irreducible.

Suppose that \(m=1\).
Set \(B=\Phi(1,e^{\pi \mathbf{i}/n})\) and \(S=\Phi(e^{\pi \mathbf{i}/(2m)},\mathbf{j})\).
With respect to the splitting \(\operatorname{span}\{1,\mathbf{i}\}\oplus \operatorname{span}\{\mathbf{j},\mathbf{k}\}\cong\mathbb R^2\oplus\mathbb R^2\),
\[
B=
\operatorname{diag}\bigl(
\operatorname{Rot}(-\pi/n),
\operatorname{Rot}(\pi/n)
\bigr),
\qquad
S(v,w)=(J_0w,-J_0v).
\]
Since \(\Phi(e^{\pi \mathbf{i}/m},1)=-I_4=B^n\), we have \(G_3=\langle B,S\rangle\).
Since \(n\ge2\), \(\pi/n\not\equiv0,\pi\pmod{2\pi}\).
By Lemma~\ref{lem:special-block-switch} with \(\delta=\pi/n\), the group \(G_3\) has no invariant line or invariant \(3\)-dimensional subspace.
Its invariant \(2\)-dimensional subspaces are precisely \(W_\theta\) for \(\theta\in\mathbb R\).
\end{proof}

\subsection{\textup{Families} \texorpdfstring{\(21\), \(21'\), \(26\), \(26'\), \(26''\), \(31\), \(31'\), \(32\), and \(32'\)}{21, 21', 26, 26', 26'', 31, 31', 32, and 32'}}\label{subsec:appendix-polyhedral}

For \(P\in\{T^*,O^*,I^*\}\), set
\[
\widetilde G_{P,+}
=
\{(p,\varepsilon p)\mid p\in P,\ \varepsilon\in\{-1,1\}\},
\qquad
\widetilde G_{P,-}
=
\{(p,p)\mid p\in P\}.
\]
The signs \(+\) and \(-\) correspond to the unprimed and primed families, respectively.
For \(P=T^*,O^*,I^*\), the pairs \(\bigl(\Phi(\widetilde G_{P,+}),\Phi(\widetilde G_{P,-})\bigr)\) are the chosen representatives for \((21,21')\), \((26,26')\), and \((31,31')\), respectively, as in Table~\ref{tab:mecchia-seppi-families-20-34}.

\begin{proposition}\label{prop:case-family-21-21p}
Let \(G<\SO(4)\) be the chosen representative in \textup{Family} \(21\) or \(21'\).
The splitting
\[
\mathbb H=\operatorname{span}\{1\}\oplus W
\]
is \(G\)-invariant.
The real \(G\)-module \(W\) is irreducible.
\end{proposition}

\begin{proof}
The chosen representatives in \(S^3\times S^3\) satisfy \(\widetilde G_{21}=\widetilde G_{T^*,+}\) and \(\widetilde G_{21'}=\widetilde G_{T^*,-}\).
By Lemma~\ref{lem:diagonal-polyhedral-module} with \(P=T^*\), the splitting \(\mathbb H=\operatorname{span}\{1\}\oplus W\) is invariant under \(G_{21}\) and \(G_{21'}\).
The real module \(W\) is irreducible for both groups.
\end{proof}

\begin{proposition}\label{prop:case-family-26}
Let \(G<\SO(4)\) be the chosen representative in \textup{Family} \(26\) or \(26'\).
The splitting
\[
\mathbb H=\operatorname{span}\{1\}\oplus W
\]
is \(G\)-invariant.
The real \(G\)-module \(W\) is irreducible.

For \textup{Family} \(26''\), the \(G_{26''}\)-module \(\mathbb H\) admits the invariant splitting
\[
\mathbb H
=
\operatorname{span}\{\mathbf{i}\}
\oplus
\operatorname{span}\{1,\mathbf{j},\mathbf{k}\}.
\]
The real \(G_{26''}\)-module \(\operatorname{span}\{1,\mathbf{j},\mathbf{k}\}\) is irreducible.
\end{proposition}

\begin{proof}
The chosen representatives in \(S^3\times S^3\) satisfy \(\widetilde G_{26}=\widetilde G_{O^*,+}\) and \(\widetilde G_{26'}=\widetilde G_{O^*,-}\).
By Lemma~\ref{lem:diagonal-polyhedral-module} with \(P=O^*\), the splitting \(\mathbb H=\operatorname{span}\{1\}\oplus W\) is invariant under \(G_{26}\) and \(G_{26'}\).
The real module \(W\) is irreducible for both groups.

For \(q\in O^*\) and \(v\in\mathbb H\), the identity \(f_0(q)=\epsilon_{26}(q)q\) implies \(\Phi(q,f_0(q))(v)=\epsilon_{26}(q)qvq^{-1}\).
Write \(\operatorname{graph}(f_0)=\{(q,f_0(q))\mid q\in O^*\}\).
The splitting \(\operatorname{span}\{1\}\oplus W\) is invariant under \(\Phi(\operatorname{graph}(f_0))\).
By Lemma~\ref{lem:diagonal-polyhedral-module}, the subgroup \(H_T<\Phi(\operatorname{graph}(f_0))\) acts irreducibly on \(W\).

Conjugation by \((1,\mathbf{i}^{-1})\) maps \(\operatorname{graph}(f_0)\) onto \(\widetilde G_{26''}\).
The map \(\Phi(1,\mathbf{i}^{-1})\colon v\mapsto v\mathbf{i}\) carries the splitting \(\operatorname{span}\{1\}\oplus W\) to
\[
\operatorname{span}\{\mathbf{i}\}
\oplus
\operatorname{span}\{1,\mathbf{j},\mathbf{k}\}.
\]
The real \(G_{26''}\)-module \(\operatorname{span}\{1,\mathbf{j},\mathbf{k}\}\) is irreducible.
\end{proof}

\begin{proposition}\label{prop:case-family-31-31p}
Let \(G<\SO(4)\) be the chosen representative in \textup{Family} \(31\) or \(31'\).
The splitting
\[
\mathbb H=\operatorname{span}\{1\}\oplus W
\]
is \(G\)-invariant.
The real \(G\)-module \(W\) is irreducible.
\end{proposition}

\begin{proof}
The chosen representatives in \(S^3\times S^3\) satisfy \(\widetilde G_{31}=\widetilde G_{I^*,+}\) and \(\widetilde G_{31'}=\widetilde G_{I^*,-}\).
By Lemma~\ref{lem:diagonal-polyhedral-module} with \(P=I^*\), the splitting \(\mathbb H=\operatorname{span}\{1\}\oplus W\) is invariant under \(G_{31}\) and \(G_{31'}\).
The real module \(W\) is irreducible for both groups.
\end{proof}

Let \(K=\mathbb Q(\sqrt5)\).
Set
\[
\mathbb H_K
=
\{a+b\mathbf{i}+c\mathbf{j}+d\mathbf{k}\mid a,b,c,d\in K\}.
\]
Let \(a\mapsto a^\dagger\) denote the nontrivial automorphism of \(K\).
Extend this automorphism coefficientwise to a \(\mathbb Q\)-algebra automorphism of \(\mathbb H_K\).
Du Val realizes \(I^*\) as a subgroup of \(\mathbb H_K\)~\cite[pp.~55--57]{val1964homographies}.
Set \((I^*)^\dagger=\{r^\dagger\mid r\in I^*\}\).
For the chosen representatives in \textup{Families} \(32\) and \(32'\), we have \(L=(I^*)^\dagger\) and \(R=I^*\).
Set
\[
\widetilde G_{32'}
=
\{(r^\dagger,r)\mid r\in I^*\},
\qquad
\widetilde G_{32}
=
\{(\varepsilon r^\dagger,r)\mid r\in I^*,\ \varepsilon\in\{\pm1\}\}.
\]

\begin{proposition}\label{prop:case-family-32-32p}
Let \(G<\SO(4)\) be the chosen representative in \textup{Family} \(32\) or \(32'\).
Then the real \(G\)-module \(\mathbb H\) is irreducible.
\end{proposition}

\begin{proof}
The coefficientwise map \(q\mapsto q^\dagger\) fixes \(\mathbf{i}\), \(\mathbf{j}\), and \(w_1\).
It therefore fixes
\[
T^*=\langle \mathbf{i},\mathbf{j},w_1\rangle
\]
pointwise.
Both \(G_{32}\) and \(G_{32'}\) contain \(H_T\).
By Lemma~\ref{lem:diagonal-polyhedral-module}, the \(H_T\)-module \(\mathbb H\) has the invariant splitting
\[
\mathbb H=\operatorname{span}\{1\}\oplus W,
\]
and the two summands are nonisomorphic irreducible real \(H_T\)-modules.

Set
\[
\tau=\frac{\sqrt5+1}{2},
\qquad
w_3=\frac12(\tau^{-1}+\tau\mathbf{j}+\mathbf{k}).
\]
Then \(w_3\in I^*\) and \(\tau^\dagger=-\tau^{-1}\).
Set
\[
\sigma_{32'}=\Phi(w_3^\dagger,w_3)\in G_{32'},
\qquad
\sigma_{32}=\Phi(w_3^\dagger,-w_3)\in G_{32}.
\]
For \(\sigma_{32}\), take \(r=-w_3\) and \(\varepsilon=-1\) in the definition of \(\widetilde G_{32}\).

Multiplication gives
\[
w_3^\dagger w_3^{-1}
=
-\frac14+\frac{\sqrt5}{4}(\mathbf{i}+\mathbf{j}+\mathbf{k})
\notin\operatorname{span}\{1\}.
\]
Thus \(\sigma_{32'}(1)=w_3^\dagger w_3^{-1}\) and \(\sigma_{32}(1)=-w_3^\dagger w_3^{-1}\).
Neither transformation preserves \(\operatorname{span}\{1\}\).
Since both transformations are orthogonal and \(W=\operatorname{span}\{1\}^\perp\), neither preserves \(W\).

By Lemma~\ref{lem:unique-one-three-splitting}, \(\operatorname{span}\{1\}\) is the unique \(H_T\)-invariant line.
The same lemma shows that \(W\) is the unique \(H_T\)-invariant \(3\)-dimensional subspace.
Suppose that \(U\) is an \(H_T\)-invariant \(2\)-dimensional subspace.
The image of \(U\) under orthogonal projection onto \(W\) is \(H_T\)-invariant.
This image has dimension at most \(2\).
Since \(W\) is irreducible, this image is zero.
This would imply \(U\subset\operatorname{span}\{1\}\), a contradiction.
Thus the only nonzero proper \(H_T\)-invariant subspaces are \(\operatorname{span}\{1\}\) and \(W\).
Every \(G_{32}\)-invariant or \(G_{32'}\)-invariant subspace is \(H_T\)-invariant.
Since \(\sigma_{32}\) and \(\sigma_{32'}\) preserve neither subspace, neither \(G_{32}\) nor \(G_{32'}\) has a nonzero proper invariant subspace.
Hence \(\mathbb H\) is irreducible as a real module for \(G_{32}\) and for \(G_{32'}\).
\end{proof}

\subsection{\textup{Families} \texorpdfstring{\(33\), \(33'\), and \(34\)}{33, 33', and 34}}\label{subsec:appendix-thirty-three-four}

\begin{remark}\label{rem:family-33p-angle-correction}
In the quotient isomorphism for \textup{Family} \(33'\) in~\cite[p.~816]{mecchia2015fibered}, the expressions \(e^{\pi\mathbf{i}/m}\) and \(e^{\pi\mathbf{i}/n}\) should be \(e^{\pi\mathbf{i}/(2m)}\) and \(e^{\pi\mathbf{i}/(2n)}\), respectively.

Set \(a=e^{\pi\mathbf{i}/(2m)}\) and \(b=e^{\pi\mathbf{i}/(2n)}\).
The coset \(a\Cyc m\) has order \(4\), whereas \(a^2\Cyc m=e^{\pi\mathbf{i}/m}\Cyc m\) has order \(2\).
For odd \(n\), the coset \(\mathbf{j}\Cyc n\) also has order \(4\).

Similarly, the coset \(b\Cyc n\) has order \(4\), whereas \(b^2\Cyc n=e^{\pi\mathbf{i}/n}\Cyc n\) has order \(2\).
For odd \(m\), the coset \(\mathbf{j}\Cyc m\) also has order \(4\).

The calculations for \textup{Family} \(33'\) use the corrected expressions.
\end{remark}

\begin{proposition}\label{prop:case-family-33-33p}
Let \(G<\SO(4)\) be the chosen representative in \textup{Family} \(33\) or \(33'\).
Let \(m,n\) be the family parameters of \(G\) as in Definition~\ref{def:representative-family-parameters}.
Then the real \(G\)-module \(\mathbb H\) is irreducible.
\end{proposition}

\begin{proof}
Set \(\kappa=1\) for \textup{Family} \(33\) and \(\kappa=2\) for \textup{Family} \(33'\).
By the description of the group in Table~\ref{tab:mecchia-seppi-families-20-34} and Proposition~\ref{prop:generator}, \(G\) contains elements acting on \(\operatorname{span}\{1,\mathbf{i}\}\oplus \operatorname{span}\{\mathbf{j},\mathbf{k}\}\) as
\begin{align*}
A&=\operatorname{diag}\bigl(\operatorname{Rot}(\kappa\pi/m),\operatorname{Rot}(\kappa\pi/m)\bigr),\\
B&=\operatorname{diag}\bigl(\operatorname{Rot}(-\kappa\pi/n),\operatorname{Rot}(\kappa\pi/n)\bigr).
\end{align*}
In both families, \(G\) contains \(S=\Phi(e^{\pi \mathbf{i}/(2m)},\mathbf{j})\), and
\[
S(\operatorname{span}\{1,\mathbf{i}\})
=
\operatorname{span}\{\mathbf{j},\mathbf{k}\}.
\]
For \textup{Family} \(33'\), the inclusion \(S\in G\) follows from the corrected quotient isomorphism in Remark~\ref{rem:family-33p-angle-correction}.

By Definition~\ref{def:representative-family-parameters}, \(m,n\ge2\) in \textup{Family} \(33\).
In \textup{Family} \(33'\), \(m,n\ge3\) are odd.
Thus
\[
\frac{\kappa\pi}{m}\not\equiv0,\pi\pmod{2\pi},
\qquad
\frac{\kappa\pi}{n}\not\equiv0,\pi\pmod{2\pi}.
\]
Lemma~\ref{lem:two-block-swap-irreducible} shows that the real \(G\)-module \(\mathbb H\) is irreducible.
\end{proof}

For \textup{Family} \(34\), set
\[
\widetilde G_{34}
=
(\Cyc{4m}/\Cyc m,\BDih{4n}/\Cyc n),
\qquad
G_{34}=\Phi(\widetilde G_{34}).
\]

\begin{proposition}\label{prop:case-family-34}
Let \(m,n\) be the family parameters of \(G_{34}\) as in Definition~\ref{def:representative-family-parameters}.
If \(m,n\ge3\), then the real \(G_{34}\)-module \(\mathbb H\) is irreducible.
If \(m=1\), then \(G_{34}\) has no invariant line or invariant \(3\)-dimensional subspace.
In this case, its invariant \(2\)-dimensional subspaces are precisely \(W_\theta\) for \(\theta\in\mathbb R\).
\end{proposition}

\begin{proof}
By Proposition~\ref{prop:generator}, \(\widetilde G_{34}\) is generated by
\[
(e^{2\pi \mathbf{i}/m},1),
\qquad
(1,e^{2\pi \mathbf{i}/n}),
\qquad
(e^{\pi \mathbf{i}/2m},\mathbf{j}).
\]
If \(m,n\ge3\), the images of these generators under \(\Phi\) satisfy the hypotheses of Lemma~\ref{lem:two-block-swap-irreducible} with \(\alpha=2\pi/m\) and \(\delta=2\pi/n\).
The real \(G_{34}\)-module \(\mathbb H\) is therefore irreducible.

Suppose that \(m=1\).
Set
\[
B=\Phi(1,e^{2\pi \mathbf{i}/n}),
\qquad
S=\Phi(e^{\pi \mathbf{i}/(2m)},\mathbf{j}).
\]
With respect to the splitting \(\operatorname{span}\{1,\mathbf{i}\}\oplus \operatorname{span}\{\mathbf{j},\mathbf{k}\}\cong\mathbb R^2\oplus\mathbb R^2\),
\[
B
=
\operatorname{diag}\bigl(
\operatorname{Rot}(-2\pi/n),
\operatorname{Rot}(2\pi/n)
\bigr),
\qquad
S(v,w)=(J_0w,-J_0v).
\]
Since \(\Phi(e^{2\pi \mathbf{i}/m},1)=I_4\), we have \(G_{34}=\langle B,S\rangle\).
The condition \(n\ge3\) gives \(2\pi/n\not\equiv0,\pi\pmod{2\pi}\).
By Lemma~\ref{lem:special-block-switch} with \(\delta=2\pi/n\), the group \(G_{34}\) has no invariant line or invariant \(3\)-dimensional subspace.
Its invariant \(2\)-dimensional subspaces are precisely \(W_\theta\) for \(\theta\in\mathbb R\).

\end{proof}

\section*{Declaration of AI use}

The second author used ChatGPT (OpenAI) for English-language editing, LaTeX assistance, reference checking, organization of the exposition, and assistance in the review and revision of selected mathematical arguments. The authors take full responsibility for all mathematical statements and proofs in this paper.

\end{document}